# Theory of vibrators with variable-order fractional forces

**Ming Li** [1, 2]

[1] Ocean College, Zhejiang University, Zhejiang 310012, P. R. China

[2] Village 1, East China Normal University, Shanghai 200062, PR. China

Emails: mli@ee.ecnu.edu.cn, ming_lihk@yahoo.com, mli15@zju.edu.cn

(Correspondence: Ocean College, Zhejiang University, China)
URL: http://orcid.org/0000-0002-2725-353X

**Abstract:** In this paper, we present a theory of six classes of vibrators with variable-order fractional forces of inertia, damping, and restoration. The novelty and contributions of the present theory are reflected in six aspects. 1) Equivalent motion equations of those variable-order fractional vibrators are proposed. 2) The analytical expressions of the effective mass, damping, and stiffness of those variable-order fractional vibrators are presented. 3) The asymptotic properties of the effective mass, damping, and stiffness of a class of variable-order fractional vibrators are given. 4) The restricted effective parameters (damping ratio, damping free natural frequency, damped natural frequency, frequency ratio) of the variable-order fractional vibrators are put forward. 5) We bring forward the analytical representations of the free responses, the impulse responses, and the frequency transfer functions of those variable-order fractional vibrators. 6) We propose a solution to an open problem of how to mathematically explain the Rayleigh damping assumption based on the present theory of variable-order fractional vibrations.

**1. Introduction**

Fractional vibrations attract the interests of researchers in various fields, ranging from physics to mathematics, see e.g., [1-20]. The literature for a class of fractional vibrators below is rich,

$$m\frac{d^\alpha x(t)}{dt^\alpha} + kx(t) = f(t), \quad 1 < \alpha \leq 2, \tag{1.1}$$

see e.g., Uchaikin [1, Chap. 7], Duan [2, Eq. (3)], Duan et al. [3, Eq. (4.2)], Al-rabtah et al. [4, Eq. (3.1)], Zurigat [5, Eq. (16)], Blaszczyk and Ciesielski [6, Eq. (1)], Blaszczyk et al. [7, Eq. (10)], Drozdov [8, Eq. (9)], Stanislavsky [9], Achar et al. ([10, Eq. (1)], [11, Eq. (9)], [12, Eq. (2)]), Tofighi [13, Eq. (2)], Ryabov and Puzenko [14, Eq. (1)], Ahmad and Elwakil [15, Eq. (1)], Blaszczyk [16], Tavazoei [17], Sandev and Tomovski [18, Eq. (36)], Singh et al. [19], Eab and Lim [20], Li et al. [21, 22], Li [23, 24].

Following the terms used in the author's recent work [23, 24], a vibrator in (1.1) is called a class I fractional one. A class I fractional vibrator is only with the fractional inertia force $mx^{(\alpha)}(t)$. Generally, $1 < \alpha \leq 2$ [1-22].

If a vibrator's motion equation takes the following form, it is called a class II fractional vibrator in [23, 24]

$$m\frac{d^2x(t)}{dt^2} + c\frac{d^\beta x(t)}{dt^\beta} + kx(t) = f(t), \quad 0 < \beta \le 1. \tag{1.2}$$

A class II fractional vibrator is simply with the fractional damping force $cx^{(\beta)}(t)$. The reports regarding the research of (1.2) are affluent, see e.g., Li [23, 24], Lin et al. [25, Eq.(2)], Duan [26, Eq. (31)], Alkhaldi et al. [27, Eq. (1a)], Dai et al. [28, Eq. (1)], Ren et al. [29, Eq. (1)], Xu et al. [30, Eq. (1)], He et al. [31, Eq. (4)], Leung et al. [32, Eq. (2)], Chen et al. [33, Eq. (1)], Deü and Matignon [34, Eq. (1)], Drăgănescu et al. [35, Eq. (4)], Rossikhin and Shitikova [36, Eq. (3)], Xie and Lin [37, Eq. (1)], Ren et al. [38, Eq. (1)], Yuan et al. [39, Eq. (8)], Dai et al. [40, Eq. (1)], Lin et al. [41, Eq. (1)], Naranjani et al. [42, Eq. (1)], Lim et al. [43], Matteo et al. [44], Tomovski and Sandev [45, Eq. (44)], Kumar et al. [46], Tian et al. [47, 48], He et al. [49, Eq. (1)], Varanis et al. [50], Pang et al. [51, Eq. (1.1)], Yu et al. [52, Eq. (2.1)], He et al. [53, Eq. (7)], Spanos and Malara [54, 55], Golmankhaneh [56], Duan et al [57, Eq. (8)].

When the motion equation of a fractional vibrator is expressed by

$$m\frac{d^\alpha x_3(t)}{dt^\alpha} + c\frac{d^\beta x_3(t)}{dt^\beta} + kx_3(t) = f(t), \quad 1 < \alpha \le 2, 0 < \beta \le 1, \tag{1.3}$$

it is called a class III fractional vibrator [23, 24]. A class III fractional vibrator is with both the fractional inertia force $mx^{(\alpha)}(t)$ and the fractional damping one $cx^{(\beta)}(t)$. For the reports about (1.3), we refer to Gomez-Aguilar et al. [58, Eq. (10)], Tian et al. [59], Berman and Cederbaum [60], Coronel-Escamilla et al. [61, Eq. (12)], Sene et al. [62, Eq. (8)], Vishwamittar et al. [63, Eq. (1)], Ismail et al. [64, Eq. (3)], Li, et al. [65].

There are other three classes of fractional vibrators that are expressed by (1.4), (1.5), and (1.6), respectively.

$$m\frac{d^\alpha x(t)}{dt^\alpha} + k\frac{d^\lambda x(t)}{dt^\lambda} = f(t), \quad 1 < \alpha < 3,\ 0 \le \lambda < 1, \tag{1.4}$$

$$m\frac{d^2 x(t)}{dt^2} + k\frac{d^\lambda x(t)}{dt^\lambda} = f(t), \quad 0 \le \lambda < 1, \tag{1.5}$$

and

$$m\frac{d^\alpha x(t)}{dt^\alpha} + c\frac{d^\beta x(t)}{dt^\beta} + k\frac{d^\lambda x(t)}{dt^\lambda} = f(t), \quad 1 < \alpha < 3,\ 0 < \beta < 2,\ 0 \le \lambda < 1. \tag{1.6}$$

The literature regarding (1.4), (1.5), and (1.6) is rare, except the recent work [24], to the best of my knowledge.

Following [24], we call (1.4) the motion equation of a class IV fractional vibrator, (1.5) for a class V fractional vibrator, and (1.6) a class VI fractional one. A class IV fractional vibrator is with both the fractional inertia force $mx^{(\alpha)}(t)$ and the fractional restoration force $kx^{(\lambda)}(t)$. A class V fractional vibrator is only with the fractional restoration force $kx^{(\lambda)}(t)$. A class VI fractional vibrator contains the fractional inertia force $mx^{(\alpha)}(t)$, the fractional damping force $cx^{(\beta)}(t)$, and the fractional restoration one $kx^{(\lambda)}(t)$ for $1 < \alpha < 3$, $0 < \beta < 2$, and $0 \leq \lambda < 1$, see [24].

Considering variable fractional orders $\alpha(\omega): [0,\infty) \to (1, 3)$, $\beta(\omega): [0,\infty) \to (1, 2)$, and $\lambda(\omega): [0,\infty) \to [0,1)$ for $(1.1) - (1.6)^1$, we have six classes of variable-order fractional vibrators respectively expressed by $(1.7) - (1.12)$ below.

$$m\frac{d^{\alpha(\omega)}x(t)}{dt^{\alpha(\omega)}} + kx(t) = f(t), \tag{1.7}$$

$$m\frac{d^2 x(t)}{dt^2} + c\frac{d^{\beta(\omega)}x(t)}{dt^{\beta(\omega)}} + kx(t) = f(t), \tag{1.8}$$

$$m\frac{d^{\alpha(\omega)}x(t)}{dt^{\alpha(\omega)}} + c\frac{d^{\beta(\omega)}x(t)}{dt^{\beta(\omega)}} + kx(t) = f(t), \tag{1.9}$$

$$m\frac{d^{\alpha(\omega)}x(t)}{dt^{\alpha(\omega)}} + k\frac{d^{\lambda(\omega)}x(t)}{dt^{\lambda(\omega)}} = f(t), \tag{1.10}$$

$$m\frac{d^2 x(t)}{dt^2} + k\frac{d^{\lambda(\omega)}x(t)}{dt^{\lambda(\omega)}} = f(t), \tag{1.11}$$

and

$$m\frac{d^{\alpha(\omega)}x(t)}{dt^{\alpha(\omega)}} + c\frac{d^{\beta(\omega)}x(t)}{dt^{\beta(\omega)}} + k\frac{d^{\lambda(\omega)}x(t)}{dt^{\lambda(\omega)}} = f(t). \tag{1.12}$$

In $(1.7) - (1.12)$, $m\frac{d^{\alpha(\omega)}x(t)}{dt^{\alpha(\omega)}}$, $c\frac{d^{\beta(\omega)}x(t)}{dt^{\beta(\omega)}}$, and $k\frac{d^{\lambda(\omega)}x(t)}{dt^{\lambda(\omega)}}$ designate variable-order fractional forces of inertia, damping, and restoration, respectively. The literature regarding $(1.7) - (1.12)$ is rarely seen.

Different classes of vibrators have their specific application areas. For instance, (1.7) implies a vibrator that is damping free in form but $m$ moves at $\frac{d^{\alpha(\omega)}x(t)}{dt^{\alpha(\omega)}}$ instead of $x''(t)$ and its displacement is $x(t)$.

In this paper, we aim at establishing a theory with respect to the variable-order fractional vibrators $(1.7) - (1.12)$. Note that $(1.7) - (1.11)$ are the special cases of (1.12). For example, when $c = 0$ and $\lambda(\omega) = 0$, (1.12) reduces to (1.7). Therefore, we detail the analysis of (1.12) in Sections 2-5. The results regarding $(1.7) - (1.11)$ are given in Section 7 as the consequences of the results from (1.12).

---

[1] In structrual vibrations, effective mass or damping or stiffness is frequency-dependent [87, 90, 109].

The present theory consists of a set of results in six aspects as follows. 1) The equivalent equations of (1.7) – (1.12) are proposed in Sections 2 and 7. 2) The analytical expressions of the effective mass, damping, and stiffness of (1.7) – (1.12) are presented in Sections 2 and 7. 3) The asymptotic properties of its effective mass, damping, and stiffness of (1.12) are put forward in Theorems 3.1-3.3 and Corollaries 3.2 and 3.3. 4) The analytic expressions of the restricted effective damping ratio, natural frequencies, and frequency ratio regarding (1.7) – (1.12) are brought forward in Sections 4 and 7. 5) The close form solutions of free responses, impulse ones, and frequency transfer functions of (1.7) – (1.12) are given in Sections 5 and 7. 6) A solution to the open problem of how to mathematically explaining the Rayleigh damping assumption is given in Theorem 6.1.

The rest of the paper is organized as follows. In Section 2, we propose an equivalent equation of (1.12) and the analytic expressions of the effective mass, damping, and stiffness of (1.12). In Section 3, we present the asymptotic properties of the effective mass, damping, and stiffness of (1.12). In Section 4, we bring forward the analytic expressions of the restricted effective damping ratio, restricted effective natural frequencies, and restricted effective frequency ratio of (1.12). In Section 5, the expressions of the free response, impulse one, and frequency transfer function regarding (1.12) are presented. In Section 6, we propose a mathematical explanation of the Rayleigh damping assumption. The corresponding results of (1.7) – (1.11) are given in Section 7, which is followed by conclusions.

## 2. Equivalent motion equation and effective mass, damping, and stiffness of fractional vibrator (1.12)

In this section, we first brief the preliminaries. Then, we present the equivalent motion equation of (1.12). Finally, we propose the analytical expressions of the effective mass, damping, and stiffness of (1.12).

### 2.1. Preliminaries

In this research, we use the Weyl fractional derivative [66-83]. Let $X(\omega)$ be the Fourier transform of $x(t)$. Using the Weyl fractional derivative, the Fourier transform of $x^{(v)}(t)$ for $v \geq 0$ is given by

$$\mathrm{F}\left[x^{(v)}(t)\right] = (i\omega)^v X(\omega), \qquad (2.1)$$

where F the Fourier transform operator, see Uchaikin [1], Miller and Ross [82], Lavoie et al, [83].

**Lemma 2.1.** Let $F_1(\omega) = \mathrm{F}[f_1(t)]$ and $F_2(\omega) = \mathrm{F}[f_2(t)]$. If $\mathrm{F}[f_1(t) - f_2(t)] = 0$, $f_1(t) - f_2(t)$ is a null function. If $f_1(t) - f_2(t)$ is a null function, $f_1(t) = f_2(t)$ in the sense of $F_1(\omega) = F_2(\omega)$ (Papoulis [84], Gelfand and Vilenkin [85], Bracewell [86]). □

### 2.2. Equivalent motion equation of (1.12)

The theorem proposed below gives an equivalent motion equation of (1.12).

**Theorem 2.1.** The motion equation of the fractional vibrator (1.12) is equivalently given by

$$-\left[m\omega^{\alpha(\omega)-2}\cos\frac{\alpha(\omega)\pi}{2}+c\omega^{\beta(\omega)-2}\cos\frac{\beta(\omega)\pi}{2}\right]\frac{d^2x(t)}{dt^2}$$
$$+\left[m\omega^{\alpha(\omega)-1}\sin\frac{\alpha(\omega)\pi}{2}+c\omega^{\beta(\omega)-1}\sin\frac{\beta(\omega)\pi}{2}+k\omega^{\lambda(\omega)-1}\sin\frac{\lambda(\omega)\pi}{2}\right]\frac{dx(t)}{dt} \quad (2.2)$$
$$+k\omega^{\lambda(\omega)}\cos\frac{\lambda(\omega)\pi}{2}x(t)=f(t),\ 1<\alpha(\omega)<3,\ 0<\beta(\omega)<2,\ 0\leq\lambda(\omega)<1.$$

*Proof.* Let $X(\omega)$ be the Fourier transform of $x(t)$. Rewrite (1.12) by

$$f_1(t) = m\frac{d^{\alpha(\omega)}x(t)}{dt^{\alpha(\omega)}} + c\frac{d^{\beta(\omega)}x(t)}{dt^{\beta(\omega)}} + k\frac{d^{\lambda(\omega)}x(t)}{dt^{\lambda(\omega)}}. \quad (2.3)$$

On one hand, (2.2) is rewritten by

$$f_2(t) = -\left[\theta m\omega^{\alpha(\omega)-2}\cos\frac{\alpha(\omega)\pi}{2}+c\omega^{\beta(\omega)-2}\cos\frac{\beta(\omega)\pi}{2}\right]\frac{d^2x(t)}{dt^2}$$
$$+\left[m\omega^{\alpha(\omega)-1}\sin\frac{\alpha(\omega)\pi}{2}+c\omega^{\beta(\omega)-1}\sin\frac{\beta(\omega)\pi}{2}+k\omega^{\lambda(\omega)-1}\sin\frac{\lambda(\omega)\pi}{2}\right]\frac{dx(t)}{dt}+k\omega^{\lambda(\omega)}\cos\frac{\lambda(\omega)\pi}{2}x(t). \quad (2.4)$$

Doing $F[f_1(t)]$ yields

$$F[f_1(t)] = \left[m(i\omega)^{\alpha(\omega)}+c(i\omega)^{\beta(\omega)}+k(i\omega)^{\lambda(\omega)}\right]X(\omega). \quad (2.5)$$

On the other side, doing $F[f_2(t)]$ produces

$$F[f_2(t)] = \left[m\omega^{\alpha(\omega)}\cos\frac{\alpha(\omega)\pi}{2}+c\omega^{\beta(\omega)}\cos\frac{\beta(\omega)\pi}{2}\right]X(\omega)$$
$$+i\left[m\omega^{\alpha(\omega)}\sin\frac{\alpha(\omega)\pi}{2}+c\omega^{\beta(\omega)}\sin\frac{\beta(\omega)\pi}{2}+k\omega^{\lambda(\omega)}\sin\frac{\lambda(\omega)\pi}{2}\right]X(\omega)+k\omega^{\lambda(\omega)}\cos\frac{\lambda(\omega)\pi}{2}X(\omega), \quad (2.6)$$

where $i = \sqrt{-1}$. Rewriting the above yields

$$F[f_2(t)] = \left[m\omega^{\alpha(\omega)}\cos\frac{\alpha(\omega)\pi}{2}+im\omega^{\alpha(\omega)}\sin\frac{\alpha(\omega)\pi}{2}\right]X(\omega)$$
$$+\left[c\omega^{\beta(\omega)}\cos\frac{\beta(\omega)\pi}{2}+ic\omega^{\beta(\omega)}\sin\frac{\beta(\omega)\pi}{2}\right]X(\omega)+\left[k\omega^{\lambda(\omega)}\cos\frac{\lambda(\omega)\pi}{2}+ik\omega^{\lambda(\omega)}\sin\frac{\lambda(\omega)\pi}{2}\right]X(\omega). \quad (2.7)$$

Since the principal values of $i^{\alpha(\omega)}$, $i^{\beta(\omega)}$, and $i^{\lambda(\omega)}$ are given by

$$i^{\alpha(\omega)} = \cos\frac{\alpha(\omega)\pi}{2}+i\sin\frac{\alpha(\omega)\pi}{2},$$
$$i^{\beta(\omega)} = \cos\frac{\beta(\omega)\pi}{2}+i\sin\frac{\beta(\omega)\pi}{2}, \quad (2.8)$$
$$i^{\lambda(\omega)} = \cos\frac{\lambda(\omega)\pi}{2}+i\sin\frac{\lambda(\omega)\pi}{2},$$

we can write (2.7) by

$$F[f_2(t)] = \left[ m(i\omega)^{\alpha(\omega)} + c(i\omega)^{\beta(\omega)} + k(i\omega)^{\lambda(\omega)} \right] X(\omega). \tag{2.9}$$

From (2.5) and (2.9), we see that $F[f_1(t)] = F[f_2(t)]$. Thus, $F[f_1(t) - f_2(t)] = 0$. Hence, $f_1(t) = f_2(t)$ as $f_1(t) - f_2(t)$ is a null function. □

**2.3. Effective mass, damping, and stiffness of (1.12)**

Note that, in (1.12), the unit of $mx^{(\alpha)}(t)$ is not Newton. Thus, the primary mass $m$ is not an effective quantity to be an inertia measure unless $\alpha = 2$. Also, the primary damping coefficient $c$ does not effectively quantify a damping measure in (1.12) since the unit of $cx^{(\beta)}(t)$ is not Newton if $\beta \neq 1$. Similarly, the primary stiffness $k$ in (1.12) is not an effective quantity to be a restoration measure as the unit of $kx^{(\lambda)}(t)$ is non-Newton when $\lambda \neq 0$. In vibrations, the effective mass of a vibrator is the coefficient of $x''(t)$, the effective damping is the coefficient of $x'(t)$, and the effective stiffness is that of $x(t)$ (Harris [87], Den Hartog [88], Timoshenko [89], Palley et al. [90], Nakagawa and Ringo [91]). Now, we present the effective mass, damping, and stiffness of (1.12) by the theorems below.

**Theorem 2.2.** Let $m_{\text{eff}}$ be the effective mass of the fractional vibrator (1.12). Then,

$$m_{\text{eff}} = -\left[ m\omega^{\alpha(\omega)-2} \cos\frac{\alpha(\omega)\pi}{2} + c\omega^{\beta(\omega)-2} \cos\frac{\beta(\omega)\pi}{2} \right]. \tag{2.10}$$

*Proof.* The above is the coefficient of $x''(t)$ in (2.2). Thus, $m_{\text{eff}}$ is an inertia measure of the vibrator (2.2) and equivalently (1.12). □

Though the unit of $m_{\text{eff}}$ is not kg unless $\alpha(\omega) = 2$ and $\beta(\omega) = 1$, it measures the inertia of a vibrator (2.2) and equivalently (1.12). The quantity $m_{\text{eff}}$ reduces to the primary $m$ when $\alpha(\omega) = 2$ and $\beta(\omega) = 1$.

**Theorem 2.3.** Denote by $c_{\text{eff}}$ the effective damping coefficient for the fractional vibrator (1.12). Then,

$$c_{\text{eff}} = m\omega^{\alpha(\omega)-1} \sin\frac{\alpha(\omega)\pi}{2} + c\omega^{\beta(\omega)-1} \sin\frac{\beta(\omega)\pi}{2} + k\omega^{\lambda(\omega)-1} \sin\frac{\lambda(\omega)\pi}{2}. \tag{2.11}$$

*Proof.* The above is the coefficient of $x'(t)$ in (2.2) that is equivalently to (1.12). □

Although the unit of $c_{\text{eff}}$ is not that of the standard damping, that is, $N \times m^{-1} \times s$, it measures the damping coefficient of a vibrator (2.2) and equivalently (1.12). The quantity $c_{\text{eff}}$ degenerates to $c$ if $\alpha(\omega) = 2$, $\beta(\omega) = 1$, and $\lambda(\omega) = 0$.

**Theorem 2.4.** Let $k_{\text{eff}}$ be the effective stiffness of the fractional vibrator (1.12). Then,

$$k_{\text{eff}} = k\omega^{\lambda(\omega)} \cos\frac{\lambda(\omega)\pi}{2}. \tag{2.12}$$

*Proof.* In (2.2), the above is the coefficient of $x(t)$ in (2.2) which is equivalently to (1.12). □

The unit of $k_{\text{eff}}$ is not $N \times m^{-1}$, but it measures the restoration of the vibrator (2.2) and equivalently (1.12). The quantity $k_{\text{eff}}$ reduces to the primary $k$ if $\lambda(\omega) = 0$.

## 3. Asymptotic properties of effective mass, damping, and stiffness of variable-order fractional vibrator (1.12)

As $m_{\text{eff}}$, $c_{\text{eff}}$, and $k_{\text{eff}}$ are the functions of vibration frequency $\omega$, we may write them by $m_{\text{eff}}(\omega)$, $c_{\text{eff}}(\omega)$, and $k_{\text{eff}}(\omega)$.

**Theorem 3.1 (Asymptotic property for $\omega \to \infty$).** If $2 < \alpha(\omega) < 3$, $m_{\text{eff}}(\omega) \to \infty$ when $\omega \to \infty$. On the other side, $m_{\text{eff}}(\omega) \to 0$ when $1 < \alpha(\omega) < 2$ if $\omega \to \infty$. That is,

$$\lim_{\omega \to \infty} m_{\text{eff}}(\omega) = \begin{cases} \infty, & 2 < \alpha(\omega) < 3 \\ 0, & 1 < \alpha(\omega) < 2 \end{cases}. \tag{3.1}$$

*Proof.* Because $\lim_{\omega \to \infty} \omega^{\alpha(\omega)-2} = \infty$ and $\cos \frac{\alpha(\omega)\pi}{2} < 0$ if $2 < \alpha(\omega) < 3$. On the other hand, $\lim_{\omega \to \infty} \omega^{\alpha(\omega)-2} = 0$ for $1 < \alpha(\omega) < 2$. Besides, $\lim_{\omega \to \infty} \omega^{\beta(\omega)-2} = 0$ for $0 < \beta(\omega) < 2$. Thus, the above is valid. □

**Theorem 3.2 (Negative mass).** If $0 < \beta(\omega) < 1$ and $2 < \alpha(\omega) < 3$, $m_{\text{eff}} \to -\infty$ if $\omega \to 0$.

*Proof.* Note that $\cos \frac{\beta(\omega)\pi}{2} > 0$ and $\lim_{\omega \to 0} \omega^{\beta-2} = \infty$ for $0 < \beta(\omega) < 1$. Besides, $\lim_{\omega \to 0} \omega^{\alpha(\omega)-2} = 0$ for $2 < \alpha(\omega) < 3$. Thus, for $0 < \beta(\omega) < 1$ and $2 < \alpha(\omega) < 3$, we have

$$\lim_{\omega \to 0} m_{\text{eff}}(\omega) = -\infty. \tag{3.2}$$

The proof finishes[2]. □

The above implies that the range of $m_{\text{eff}}(\omega)$ is $(-\infty, \infty)$ in general. A few of plots of $m_{\text{eff}}(\omega)$ are shown in Fig. 3.1.

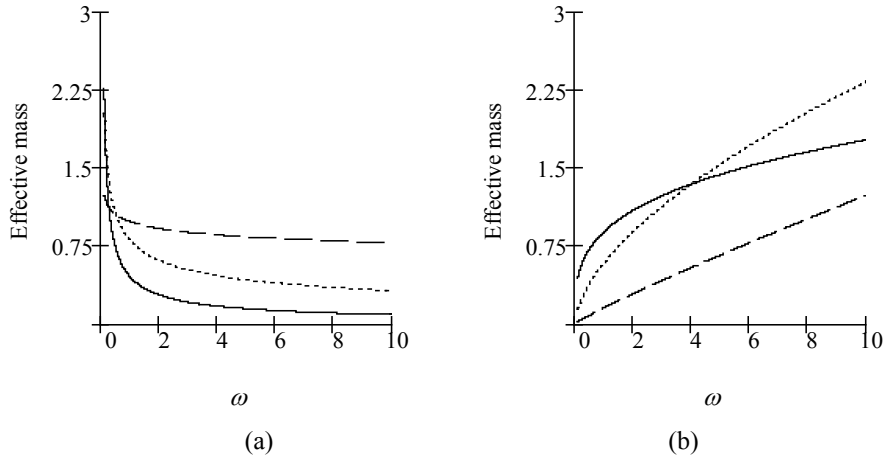

(a)            (b)

---

[2] Refer [24, Chap. 15] and references therein about negative mass.

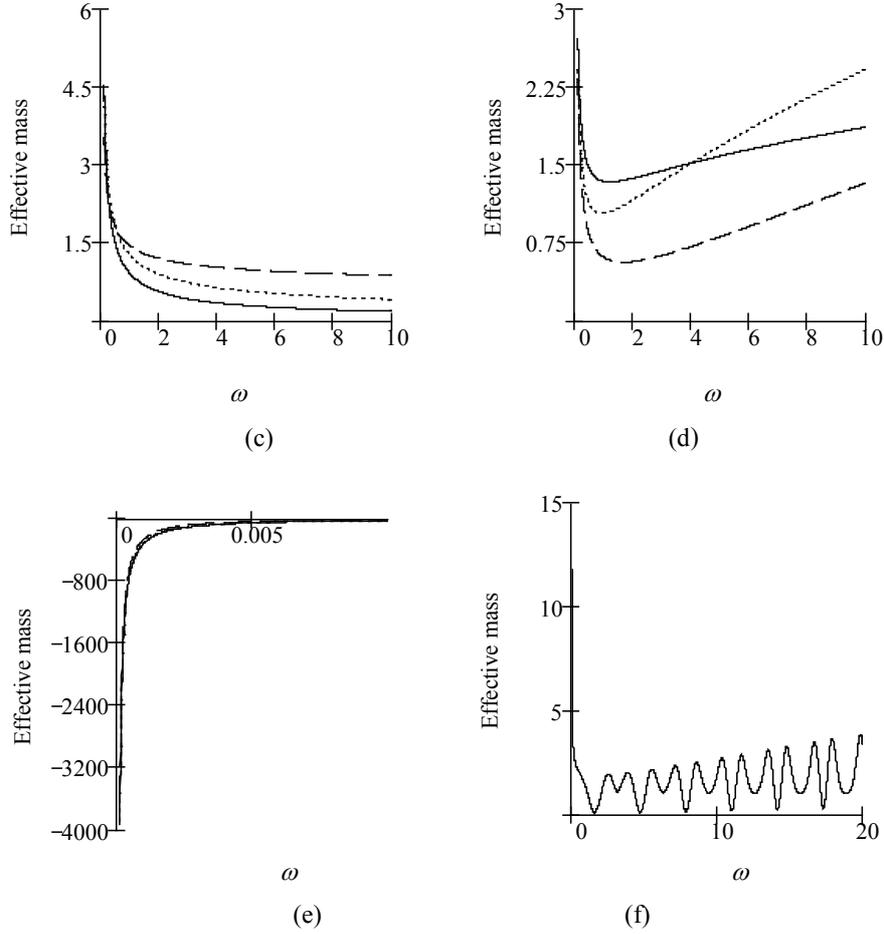

Fig. 3.1. Plots of $m_{\text{eff}}(\omega)$ with $m = 1$, $c = 1$, and $k = 1$. (a). $\beta = 1$. Solid: $\alpha = 1.3$. Dot: $\alpha = 1.6$. Dash: $\alpha = 1.9$. (b). $\beta = 1$. Solid: $\alpha = 2.3$. Dot: $\alpha = 2.6$. Dash: $\alpha = 2.9$. (c). $\beta = 1.3$. Solid: $\alpha = 1.3$. Dot: $\alpha = 1.6$. Dash: $\alpha = 1.9$. (d). $\beta = 1.3$. Solid: $\alpha = 2.3$. Dot: $\alpha = 2.6$. Dash: $\alpha = 2.9$. (e). Observing negative $m_{\text{eff}}$ at small $\omega$ when $\beta = 0.9$. Solid: $\alpha = 1.3$. Dot: $\alpha = 1.6$. Dash: $\alpha = 1.9$. (f). $\alpha(\omega) = 1.10 + 1.89|\sin(\omega)|$, $\beta(\omega) = 1 + 0.99|\cos(\omega)|$.

We now consider the asymptotic property of $c_{\text{eff}}$ for $\omega \to \infty$.

**Theorem 3.3.** If $\omega \to \infty$,

$$c_{\text{eff}} = \begin{cases} \infty, 1 < \alpha(\omega) < 2, 0 < \beta(\omega) < 2 \\ -\infty, 2 < \alpha(\omega) < 3, 0 < \beta(\omega) < 1 \end{cases}. \tag{3.3}$$

*Proof.* Since $0 \le \lambda(\omega) < 1$, we have $\lim_{\omega \to \infty} \omega^{\lambda(\omega)-1} = 0$. In addition, because $\sin\frac{\beta(\omega)\pi}{2} > 0$ when $0 < \beta(\omega) < 2$ and $\sin\frac{\alpha(\omega)\pi}{2} > 0$ if $1 < \alpha(\omega) < 2$, $c_{\text{eff}}(\omega) \to \infty$ if $\omega \to \infty$ for $1 < \alpha(\omega) < 2$ and $0 < \beta(\omega) < 2$. Additionally, because $\lim\limits_{\omega\to\infty}\omega^{\beta(\omega)-1}=0$ when $0<\beta(\omega)<1$ but $\sin\dfrac{\alpha(\omega)\pi}{2}<0$ if $2<\alpha(\omega)<3$, we have $c_{\text{eff}}(\omega)\to-\infty$ for $2<\alpha(\omega)<3$ and $0<\beta(\omega)<1$ when $\omega\to\infty$. □

The above exhibits that $-\infty<c_{\text{eff}}<\infty^3$.

**Corollary 3.1.** The fractional vibrator (1.12) may be self-vibrated and accordingly non-stable if $2<\alpha(\omega)<3$ and $0<\beta(\omega)<1$.

*Proof.* As $c_{\text{eff}}$ may be negative if $2<\alpha(\omega)<3$ and $0<\beta(\omega)<1$, a vibrator (1.12) may be self-vibrated and non-stable in that case. □

**Corollary 3.2.** If $1<\alpha(\omega)<3$, $0<\beta(\omega)<1$ and $\omega$ is large enough, we have

$$c_{\text{eff}} \cong m\omega^{\alpha(\omega)-1}\sin\dfrac{\alpha(\omega)\pi}{2}. \tag{3.4}$$

In addition, for $1<\alpha(\omega)<3$ and $0<\beta(\omega)<1$, if $\omega$ is small enough,

$$c_{\text{eff}} \cong c\omega^{\beta(\omega)-1}\sin\dfrac{\beta(\omega)\pi}{2}+k\omega^{\lambda(\omega)-1}\sin\dfrac{\lambda(\omega)\pi}{2}. \tag{3.5}$$

*Proof.* If $1<\alpha(\omega)<3$, $0<\beta(\omega)<1$ and $\omega$ is sufficiently large, $\omega^{\beta(\omega)-1}\approx 0$ and $\omega^{\lambda(\omega)-1}\approx 0$. Thus, (3.4) holds. On the other hand, if $\omega$ is small enough, $\omega^{\alpha(\omega)-1}\approx 0$. Thus, (3.5) results. □

Fig. 3.2 gives a few illustrations of $c_{\text{eff}}$.

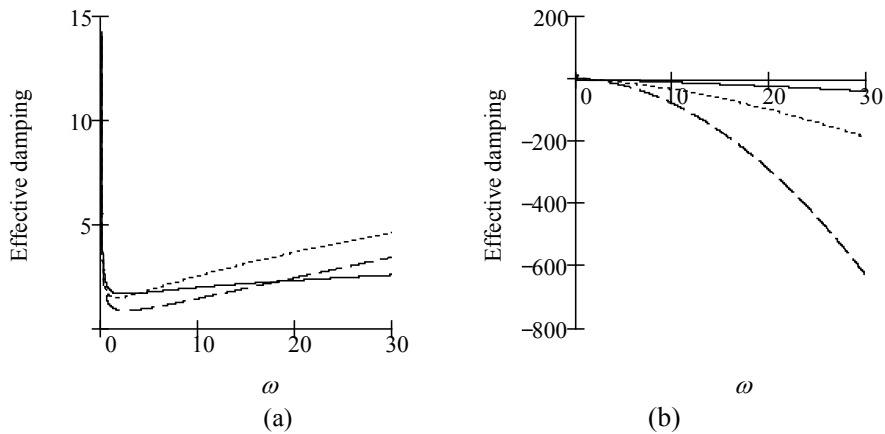

(a)    (b)

---

[3] For several particular cases of negative damping in engineering, refer Den Harton [88] and Nakagawa and Ringo [91]. For general cases of negative damping in fractional vibrations with constant fractional orders, refer to [24].

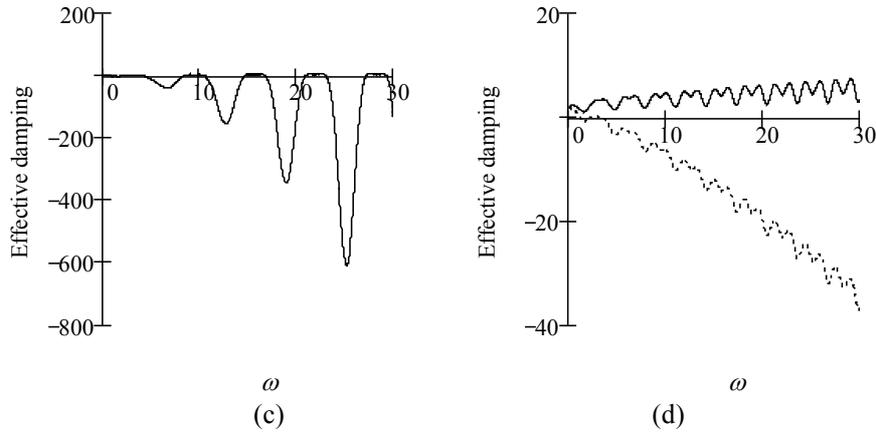

Fig. 3.2. Plots of $c_{\text{eff}}$ with $m = 1$, $c = 1$, and $k = 1$. (a). $\beta = 0.3$, $\lambda = 0.3$. Solid: $\alpha = 1.3$. Dot: $\alpha = 1.6$. Dash: $\alpha = 1.9$. (b). Negative damping with constant fractional orders, $\beta = 0.3$, $\lambda = 0.3$. Solid: $\alpha = 2.3$. Dot: $\alpha = 2.6$. Dash: $\alpha = 2.9$. (c). Negative damping with variable fractional orders, $\alpha(\omega) = 1.10 + 1.89|\cos(0.5\omega)|$, $\beta(\omega) = 1 + 0.99|\sin(\omega)|$, $\lambda(\omega) = 0.99|\sin(\omega)|$. (d). $\beta(\omega) = 1 + 0.99|\sin(\omega)|$, $\lambda(\omega) = 0.99|\cos(\omega)|$. Solid: $\alpha = 1.3$. Dot (negative damping): $\alpha = 2.3$.

**Corollary 3.3.** For $k_{\text{eff}}$, we have $k_{\text{eff}}(\omega) \geq 0$. Besides, $\lim\limits_{\omega \to \infty} \omega^{\lambda(\omega)} = \infty$ and $\lim\limits_{\omega \to 0} \omega^{\lambda(\omega)} = 0$.

*Proof.* The proof is straightforward as $0 \leq \lambda(\omega) < 1$ in (2.12). □

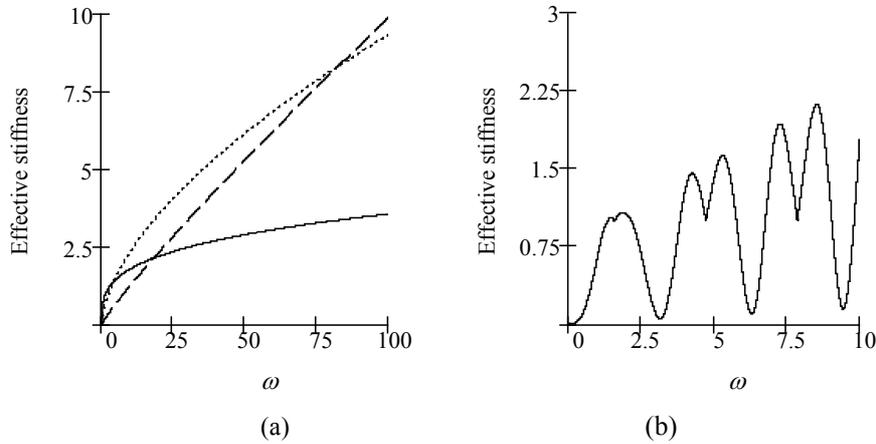

Fig. 3.3. Plots of $k_{\text{eff}}$ with $k = 1$. (a). Solid: $\lambda = 0.3$. Dot: $\lambda = 0.6$. Dash: $\lambda = 0.9$. (b). $\lambda(\omega) = 0.99|\cos(\omega)|$.

We show its plots in Fig. 3.3.

## 4. Restricted effective damping ratio, natural frequencies, and frequency ratio of variable-order fractional vibrator (1.12)

The effective parameters $m_{\text{eff}}$, $c_{\text{eff}}$, and $k_{\text{eff}}$ are directly obtained from (2.2). In this section, we address the restricted effective damping ratio, natural frequencies, and frequency ratio. By restricted, we mean that additional conditions are considered for expressing those parameters.

### 4.1. Restricted effective damping ratio of (1.12)

As can be seen from Section 3, $m_{\text{eff}} \in (-\infty, \infty)$, $c_{\text{eff}} \in (-\infty, \infty)$, and $k_{\text{eff}} \in (0, \infty)$. Let $\zeta_{\text{eff}}$ be the restricted effective damping ratio (effective damping ratio for short) of the fractional vibrator (1.12). Here, we only consider $\zeta_{\text{eff}}$ in the case of $m_{\text{eff}} > 0$ from a view of engineering. With that restriction, we define $\zeta_{\text{eff}}$ by

$$\zeta_{\text{eff}} = \frac{c_{\text{eff}}}{2\sqrt{m_{\text{eff}} k_{\text{eff}}}}. \tag{4.1}$$

**Theorem 4.1.** The quantity $\zeta_{\text{eff}}$ is in the form

$$\zeta_{\text{eff}} = \frac{m\omega^{\alpha(\omega)-1} \sin\frac{\alpha(\omega)\pi}{2} + c\omega^{\beta(\omega)-1} \sin\frac{\beta(\omega)\pi}{2} + k\omega^{\lambda(\omega)-1} \sin\frac{\lambda(\omega)\pi}{2}}{2\sqrt{-\left(m\omega^{\alpha(\omega)-2} \cos\frac{\alpha(\omega)\pi}{2} + c\omega^{\beta(\omega)-2} \cos\frac{\beta(\omega)\pi}{2}\right) k\omega^{\lambda(\omega)} \cos\frac{\lambda(\omega)\pi}{2}}}. \tag{4.2}$$

*Proof.* Substituting $m_{\text{eff}}$, $c_{\text{eff}}$, and $k_{\text{eff}}$ into (4.1) results in (4.2). □

Although $\zeta_{\text{eff}}$ is not dimensionless in general, it effectively reflects the damping ratio of (2.2) and equivalently (1.12) in the case of $m_{\text{eff}} > 0$. It is dimensionless if $\alpha(\omega) = 2$, $\beta(\omega) = 1$, and $\lambda(\omega) = 0$. In that case, it equals to the standard damping ratio since $\zeta_{\text{eff}}|_{\alpha(\omega)=2, \beta(\omega)=1, \lambda(\omega)=0} = \frac{c}{2\sqrt{mk}} = \zeta$.

A few of plots of $\zeta_{\text{eff}}$ are indicated in Fig. 4.1.

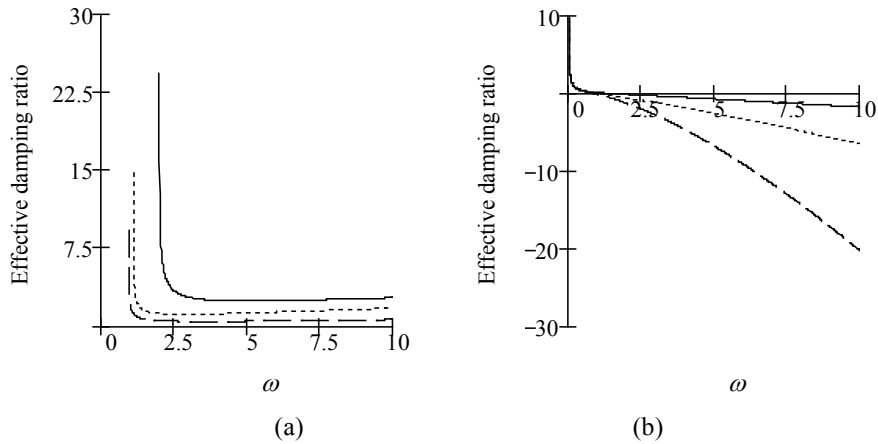

(a)      (b)

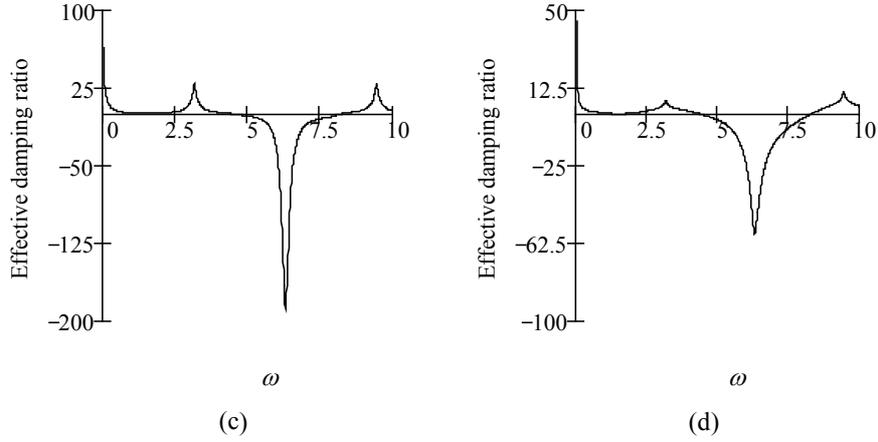

(c)             (d)

Fig. 4.1. Plots of $\zeta_{\text{eff}}$ for $m = 1 = k = c = 1$. (a). Damping ratio for some range of frequency with constant fractional orders. $\beta = 0.3$, $\lambda = 0.3$. Solid: $\alpha = 1.3$. Dot: $\alpha = 1.6$. Dash: $\alpha = 1.9$. (b). Negative damping ratio for some range of frequency with constant fractional orders. $\beta = 1.9$, $\lambda = 0.3$. Solid: $\alpha = 2.3$. Dot: $\alpha = 2.6$. Dash: $\alpha = 2.9$. (c). Negative damping ratio for some range of frequency with variable fractional orders. $\alpha(\omega) = 1.10 + 1.89|\cos(0.1\omega)|$, $\beta(\omega) = 1 + 0.99|\sin(\omega)|$, $\lambda(\omega) = 0.99|\cos(\omega)|$. (d). Negative damping ratio for some range of frequency with variable fractional orders. $\alpha(\omega) = 1.10 + 1.89|\cos(0.1\omega)|$, $\beta(\omega) = 1 + 0.99|\sin(\omega)|$, and $\lambda(\omega) = 0.99|\exp(-\omega)|$.

**4.2. Restricted effective damping free natural frequency of (1.12)**

For the fractional vibrator (1.12), with the restriction $m_{\text{eff}} > 0$ from the point of view of engineering, we coin a term restricted effective damping free natural frequency (effective damping free natural frequency in short) to the quantity defined by

$$\omega_{\text{effn}} = \sqrt{\frac{k_{\text{eff}}}{m_{\text{eff}}}}. \tag{4.3}$$

**Theorem 4.2.** The quantity $\omega_{\text{effn}}$ is in the form

$$\omega_{\text{effn}} = \sqrt{\frac{k\omega^{\lambda(\omega)}\cos\frac{\lambda(\omega)\pi}{2}}{-\left(m\omega^{\alpha(\omega)-2}\cos\frac{\alpha(\omega)\pi}{2} + c\omega^{\beta(\omega)-2}\cos\frac{\beta(\omega)\pi}{2}\right)}}. \tag{4.4}$$

*Proof.* Substituting $m_{\text{eff}}$ and $k_{\text{eff}}$ into (4.3) yields (4.4). □

Though the unit of $\omega_{\text{effn}}$ is not rad/s, but it in functional takes the form of the conventional damping free natural frequency. Its unit reduces to rad/s if $\alpha(\omega) = 2$, $\beta(\omega) = 1$, and $\lambda(\omega) = 0$ as

$$\omega_{\text{effn}}\big|_{\alpha(\omega)=2, \beta(\omega)=1, \lambda(\omega)=0} = \omega_{\text{n}} = \sqrt{\frac{k}{m}}. \tag{4.5}$$

In general, $\omega_{\text{effn}}$ is a function of $\omega$. Thus, $\omega_{\text{effn}} = \omega_{\text{effn}}(\omega)$. Fig. 4.2 illustrates a few plots of $\omega_{\text{effn}}$. For variable fractional orders, $\omega_{\text{effn}}(\omega)$ is complicated, see Figs. 4.2 (e) and (f).

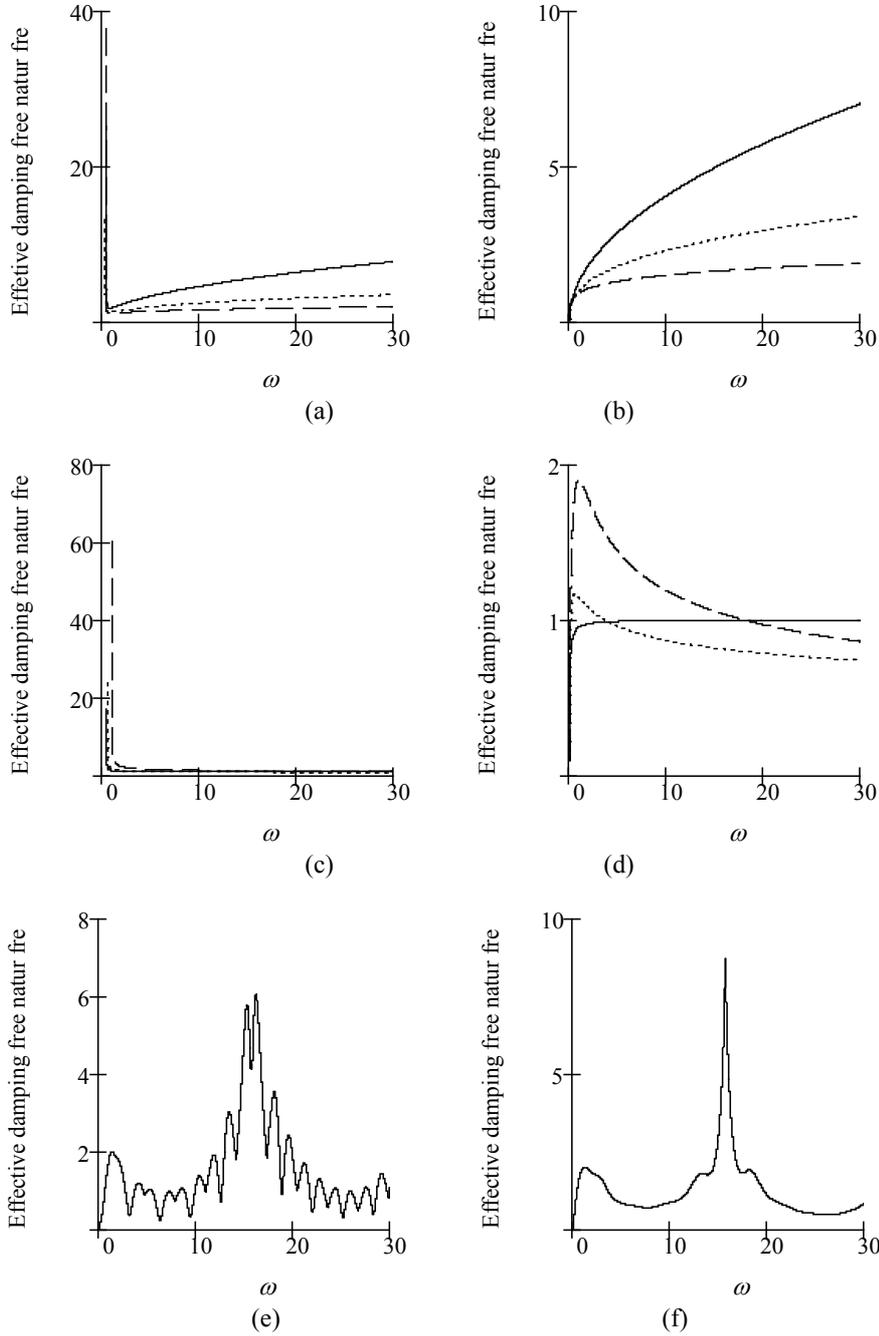

Fig. 4.2. Plots of $\omega_{\text{effn}}$ with $m = 1$ and $k = 1$, $c = 0.2$. (a). $\omega_{\text{effn}}$ with constant fractional orders. $\beta = 0.3$, $\lambda = 0.3$. Solid: $\alpha = 1.3$. Dot: $\alpha = 1.6$. Dash: $\alpha = 1.9$. (b). $\omega_{\text{effn}}$ with constant fractional orders. $\beta = 1.3$, $\lambda = 0.3$. Solid: $\alpha = 1.3$. Dot: $\alpha = 1.6$. Dash: $\alpha = 1.9$. (c). $\omega_{\text{effn}}$ with constant fractional orders. $\beta = 0.3$, $\lambda = 0.3$. Solid: $\alpha = 2.3$. Dot: $\alpha = 2.6$. Dash: $\alpha = 2.9$. (d). $\omega_{\text{effn}}$ with constant fractional orders. $\beta = 1.3$, $\lambda = 0.3$. Solid: $\alpha = 2.3$. Dot: $\alpha = 2.6$. Dash: $\alpha = 2.9$. (e). $\omega_{\text{effn}}$ with variable fractional orders. $\alpha(\omega) = 1.10 + 1.89|\cos(0.1\omega)|$,

$\beta(\omega) = 1 + 0.99|\sin(\omega)|$, $\lambda(\omega) = 0.99|\cos(\omega)|$. (f). $\omega_{\text{effn}}$ with variable fractional orders. $\alpha(\omega) = 1.10 + 1.89|\cos(0.1\omega)|$, $\beta(\omega) = 1 + 0.99|\sin(\omega)|$, $\lambda(\omega) = 0.99|\exp(-\omega)|$.

**4.3. Restricted effective damped natural frequency of (1.12)**

Now we introduce a term restricted effective damped natural frequency (effective damped natural frequency in short) for (1.12). Denote it by $\omega_{\text{effd}}$. With the restriction $|\zeta_{\text{eff}}| \leq 1$[4], we define it by

$$\omega_{\text{effd}} = \omega_{\text{effn}} \sqrt{1 - \varsigma_{\text{eff}}^2}, \quad |\varsigma_{\text{eff}}| \leq 1. \tag{4.6}$$

**Theorem 4.3.** The representation of $\omega_{\text{effd}}$ is in the form

$$\omega_{\text{effd}} = \sqrt{\frac{k\omega^{\lambda(\omega)} \cos\frac{\lambda(\omega)\pi}{2}}{-\left(m\omega^{\alpha(\omega)-2} \cos\frac{\alpha(\omega)\pi}{2} + c\omega^{\beta(\omega)-2} \cos\frac{\beta(\omega)\pi}{2}\right)}} \sqrt{1 - \left(\frac{m\omega^{\alpha(\omega)-1} \sin\frac{\alpha(\omega)\pi}{2} + c\omega^{\beta(\omega)-1} \sin\frac{\beta(\omega)\pi}{2} + k\omega^{\lambda(\omega)-1} \sin\frac{\lambda(\omega)\pi}{2}}{2\sqrt{-\left(m\omega^{\alpha(\omega)-2} \cos\frac{\alpha(\omega)\pi}{2} + c\omega^{\beta(\omega)-2} \cos\frac{\beta(\omega)\pi}{2}\right) k\omega^{\lambda(\omega)} \cos\frac{\lambda(\omega)\pi}{2}}}\right)^2}. \tag{4.7}$$

*Proof.* Substituting $\zeta_{\text{eff}}$ and $\omega_{\text{effn}}$ into (4.6) yields (4.7). □

Note that if $\alpha(\omega) = 2$, $\beta(\omega) = 1$, and $\lambda(\omega) = 0$, $\omega_{\text{effd}}$ degenerates the conventional $\omega_d$ with the unit of rad/s. It is generally a function of $\omega$. Hence, $\omega_{\text{effd}} = \omega_{\text{effd}}(\omega)$.

**4.4. Restricted effective frequency ratio of (1.12)**

Introduce a term restricted effective frequency ratio (effective frequency ratio in short) for the fractional vibrator (1.12). Denote it by $\gamma_{\text{eff}}$. Define it by

$$\gamma_{\text{eff}} = \frac{\omega}{\omega_{\text{effn}}}. \tag{4.8}$$

Both $\omega_{\text{effn}}$ and $\gamma_{\text{eff}}$ are restricted by $m_{\text{eff}} > 0$[5].

**Theorem 4.4.** The representation of $\gamma_{\text{eff}}$ is given by

$$\gamma_{\text{eff}} = \gamma \sqrt{\frac{-\left(\omega^{\alpha(\omega)-2} \cos\frac{\alpha(\omega)\pi}{2} + 2\varsigma\omega_n \omega^{\beta(\omega)-2} \cos\frac{\beta(\omega)\pi}{2}\right)}{\omega^{\lambda(\omega)} \cos\frac{\lambda(\omega)\pi}{2}}}, \tag{4.9}$$

where $\gamma = \dfrac{\omega}{\omega_n}$.

---

[4] Small damping is commonly assumed in vibration engineering, see [24, 87-92, 109]. The quantity $\omega_{\text{effd}}$ is restricted in the sense of small damping.
[5] That restriction of $m_{\text{eff}} > 0$ is taken from the point of view of vibration engineering.

*Proof.* Substituting $\omega_{\text{effn}}$ into (4.8) results in (4.9). □

The quantity $\gamma_{\text{eff}}$ reduces to $\gamma$ when $\alpha(\omega) = 2$, $\beta(\omega) = 1$, and $\lambda(\omega) = 0$. Its plots are indicated in Fig. 4.3.

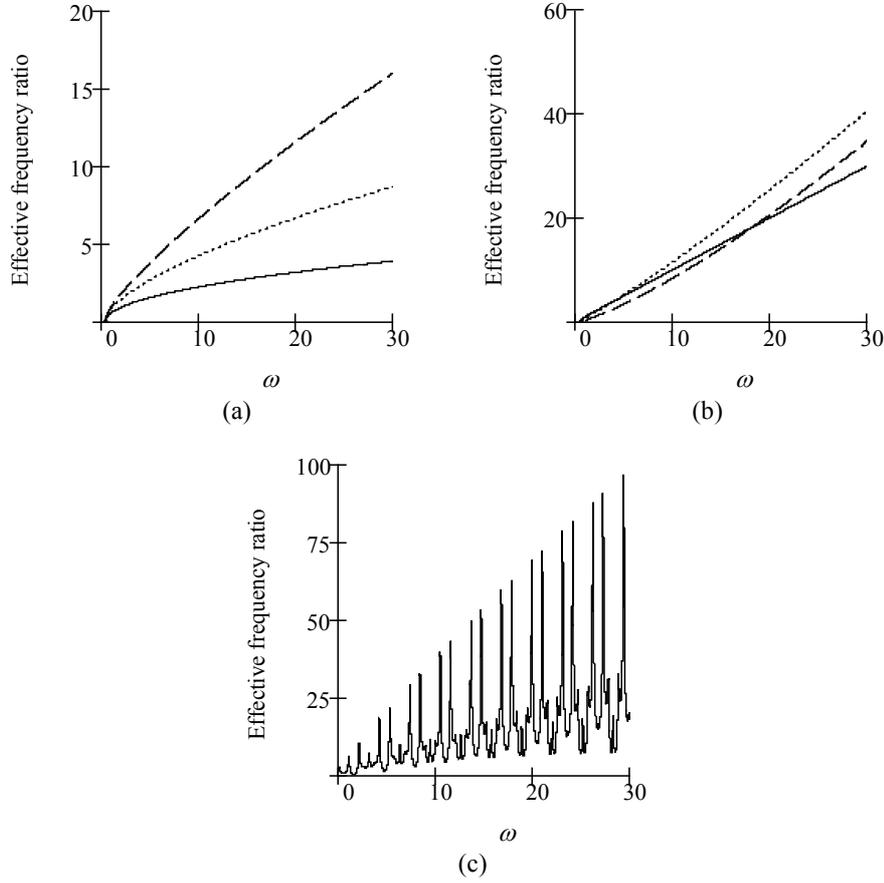

Fig. 4.3. Plots of $\gamma_{\text{eff}}$ with $m = 1$, $k = 1$, $c = 0.2$. (a). Effective frequency ratio with constant fractional orders. $\beta = 0.3$, $\lambda = 0.3$. Solid: $\alpha = 1.3$. Dot: $\alpha = 1.6$. Dash: $\alpha = 1.9$. (b). Effective frequency ratio with constant fractional orders. $\beta = 0.3$, $\lambda = 0.3$. Solid: $\alpha = 2.3$. Dot: $\alpha = 2.6$. Dash: $\alpha = 2.9$. (c). Effective frequency ratio with variable fractional orders. $\alpha(\omega) = 1.10 + 1.89|\cos(0.1\omega)|$, $\beta(\omega) = 1 + 0.99|\sin(\omega)|$, $\lambda(\omega) = 0.99|\cos(\omega)|$.

## 5. Restricted responses (free, impulse) of the variable-order fractional vibrator (1.12)

Here and below, we use two restrictions, namely, $m_{\text{eff}} > 0$ and $|\zeta_{\text{eff}}| \leq 1$, unless otherwise stated.

### 5.1. Equivalent representation of (2.2)

Based on $m_{\text{eff}}(\omega)$, $c_{\text{eff}}(\omega)$, and $k_{\text{eff}}(\omega)$, we rewrite (2.2) by

$$m_{\text{eff}} \frac{d^2 x(t)}{dt^2} + c_{\text{eff}} \frac{dx(t)}{dt} + k_{\text{eff}} x(t) = f(t). \tag{5.1}$$

The above can be further written by

$$\frac{d^2x(t)}{dt^2} + 2\varsigma_{\text{eff}}\omega_{\text{effn}}\frac{dx(t)}{dt} + \omega_{\text{effn}}^2 x(t) = \frac{f(t)}{m_{\text{eff}}}. \tag{5.2}$$

**5.2. Restricted free response**

Let $x(t)$ be the restricted free response (free response for short) of the fractional vibrator (1.12). It is a solution to the differential equation given by

$$\begin{cases} \dfrac{d^2x(t)}{dt^2} + 2\varsigma_{\text{eff}}\omega_{\text{effn}}\dfrac{dx(t)}{dt} + \omega_{\text{effn}}^2 x(t) = 0, \\ x(0) = x_0, x'(0) = v_0. \end{cases} \quad 1 < \alpha(\omega) < 3, 0 < \beta(\omega) < 2, 0 \le \lambda(\omega) < 1. \tag{5.3}$$

The theorem proposed below provides its solution.

**Theorem 5.1.** The free response $x(t)$ to (1.12) is given by

$$x(t) = e^{-\varsigma_{\text{eff}}\omega_{\text{effn}}t}\left(x_0 \cos\omega_{\text{effd}}t + \frac{v_0 + \varsigma_{\text{eff}}\omega_{\text{effn}}x_0}{\omega_{\text{effd}}}\sin\omega_{\text{effd}}t\right), \quad t \ge 0. \tag{5.4}$$

*Proof.* According to Theorem 2.1, (2.2) is equivalent to (1.12). On the other hand, (5.2) is equivalent to (2.2). Since (5.3) is a standard vibration equation in form, (5.4) holds. □

Fig. 5.1 shows a few plots of $x(t)$.

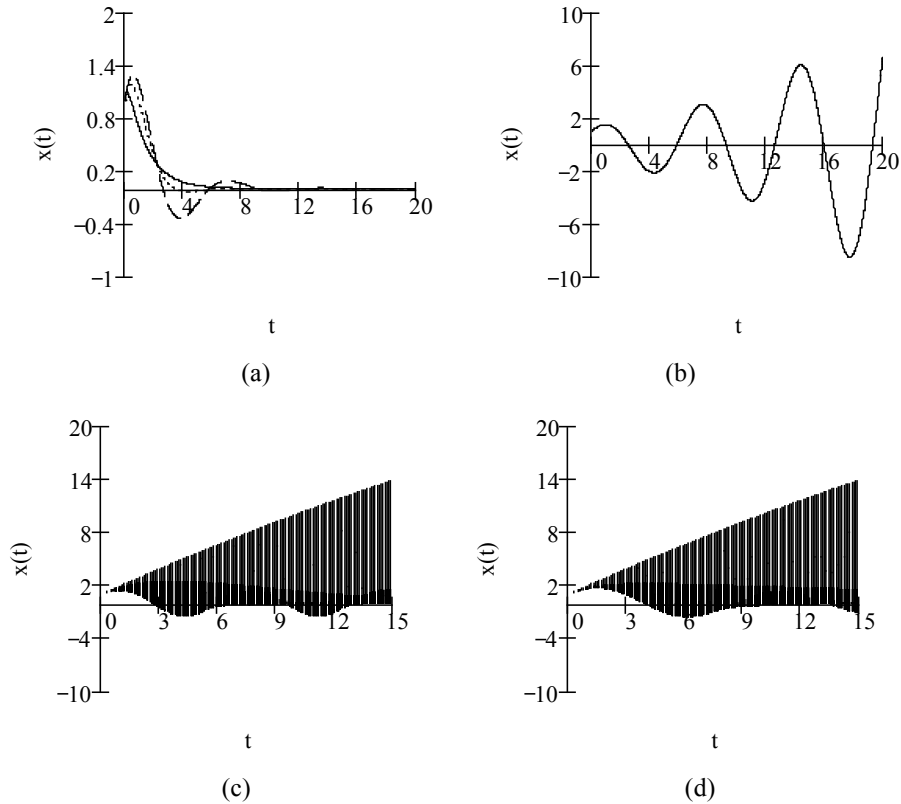

(a)  (b)  (c)  (d)

Fig. 5.1. Plots of $x(t)$ for $x_0 = 1$, $v_0 = 1$, $m = 1$ and $k = 1$. (a). For $\omega = 1.1$ and $c = 0.2$, $\beta = 0.3$, $\lambda = 0.3$. Solid: $\alpha = 1.3$. Dot. $\alpha = 1.6$. Dash: $\alpha = 1.9$. (b). Self-vibration. $x(t)$ with constant fractional orders for $\omega = 1.1$ and $c = 0.2$, $\beta = 1.9$, $\lambda = 0.3$. $\alpha = 2.4$. (c). $x(t)$ with variable fractional orders for $\omega \in (0, 1)$ and $c = 0.2$, $\alpha(\omega) = 1.10 + 1.89|\cos(0.1\omega)|$, $\beta(\omega) = 1 + 0.99|\sin(\omega)|$, $\lambda(\omega) = 0.99|\cos(\omega)|$. (d). $x(t)$ with variable fractional orders for $\omega \in (0, 1)$ and $c = 1.2$, $\alpha(\omega) = 1.10 + 1.89|\cos(0.1\omega)|$, $\beta(\omega) = 1 + 0.99|\sin(\omega)|$, $\lambda(\omega) = 0.99|\cos(\omega)|$.

In time-frequency plane, we have

$$x(t) = x(t, \omega). \tag{5.5}$$

### 5.3. Restricted impulse response

The present theorem below gives the expression of the restricted impulse response (impulse response for short) to the fractional vibrator (1.12).

**Theorem 5.2.** Let $h(t)$ be the impulse response to the fractional vibrator (1.12). It is a solution to the following fractional differential equation:

$$m\frac{d^{\alpha(\omega)}h(t)}{dt^{\alpha(\omega)}} + c\frac{d^{\beta(\omega)}h(t)}{dt^{\beta(\omega)}} + k\frac{d^{\lambda(\omega)}h(t)}{dt^{\lambda(\omega)}} = \delta(t), \quad 1 < \alpha(\omega) < 3, 0 < \beta(\omega) < 2, 0 \leq \lambda(\omega) < 1 \tag{5.6}$$

with zero initial conditions. The expression of $h(t)$ is given by

$$h(t) = e^{-\varsigma_{\text{eff}}\omega_{\text{effn}}t}\frac{1}{m_{\text{eff}}\omega_{\text{effd}}}\sin\omega_{\text{effd}}t, \quad t \geq 0. \tag{5.7}$$

*Proof.* According to Theorem 2.1, (5.6) is equivalently expressed by

$$m_{\text{eff}}\frac{d^2h(t)}{dt^2} + c_{\text{eff}}\frac{dh(t)}{dt} + k_{\text{eff}}h(t) = \delta(t).$$

The above is a standard vibration equation in form. Its solution is (5.7). □

Fig. 5.2 shows a few plots of $h(t)$.

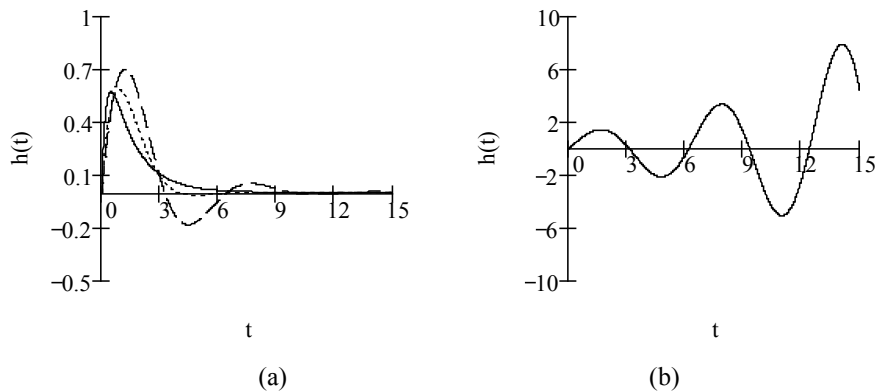

(a)     (b)

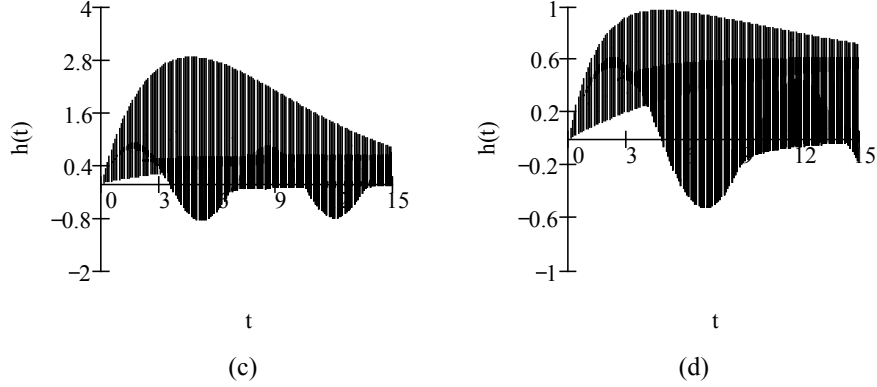

Fig. 5.2. Plots of $h(t)$ for $m = 1$ and $k = 1$. (a). $h(t)$ with constant fractional orders for $\omega = 1.1$ and $c = 0.2$, $\beta = 0.3$, $\lambda = 0.3$. Solid: $\alpha = 1.3$. Dot: $\alpha = 1.6$. Dash: $\alpha = 1.9$. (b). Self-vibration. $h(t)$ with constant fractional orders for $\omega = 1.1$ and $c = 0.2$, $\alpha = 2.5$, $\beta = 1.5$, $\lambda = 0.3$. (c). $h(t)$ with variable fractional orders for $\omega \in (0, 1)$ and $c = 0.2$, $\alpha(\omega) = 1.10 + 1.89|\cos(0.1\omega)|$, $\beta(\omega) = 1 + 0.99|\sin(\omega)|$, $\lambda(\omega) = 0.99|\cos(\omega)|$. (d). $h(t)$ with variable fractional orders for $\omega \in (0, 1)$ and $c = 1.2$, $\alpha(\omega) = 1.10 + 1.89|\cos(0.1\omega)|$, $\beta(\omega) = 1 + 0.99|\sin(\omega)|$, $\lambda(\omega) = 0.99|\cos(\omega)|$.

In time-frequency plan, $h(t)$ may be written as

$$h(t) = h(t, \omega). \tag{5.8}$$

### 5.4. Restricted frequency transfer function

**Theorem 5.3.** Let $H(\omega)$ be the restricted frequency transfer function (frequency transfer function in short) of the fractional vibrator (1.12). It is given by

$$H(\omega) = \frac{1}{k_{\text{eff}}\left(1 - \gamma_{\text{eff}}^2 + i2\varsigma_{\text{eff}}\gamma_{\text{eff}}\right)}. \tag{5.9}$$

*Proof.* Doing the Fourier transform of (5.7) produces the above. □

Using polar coordinates, we have

$$H(\omega) = |H(\omega)|e^{i\varphi(\omega)}, \tag{5.10}$$

where the amplitude-frequency response $|H(\omega)|$ is in the form

$$|H(\omega)| = \frac{1}{k_{\text{eff}}} \frac{1}{\sqrt{\left(1 - \gamma_{\text{eff}}^2\right)^2 + \left(2\varsigma_{\text{eff}}\gamma_{\text{eff}}\right)^2}}, \tag{5.11}$$

and the phase-frequency response function $\varphi(\omega)$ is expressed by

$$\varphi(\omega) = \cos^{-1} \frac{1-\gamma_{\text{eff}}^2}{\sqrt{\left(1-\gamma_{\text{eff}}^2\right)^2 + \left(2\varsigma_{\text{eff}}\gamma_{\text{eff}}\right)^2}}. \tag{5.12}$$

Figs. 5.3 and 5.4 illustrate several plots of $|H(\omega)|$ and $\varphi(\omega)$, respectively.

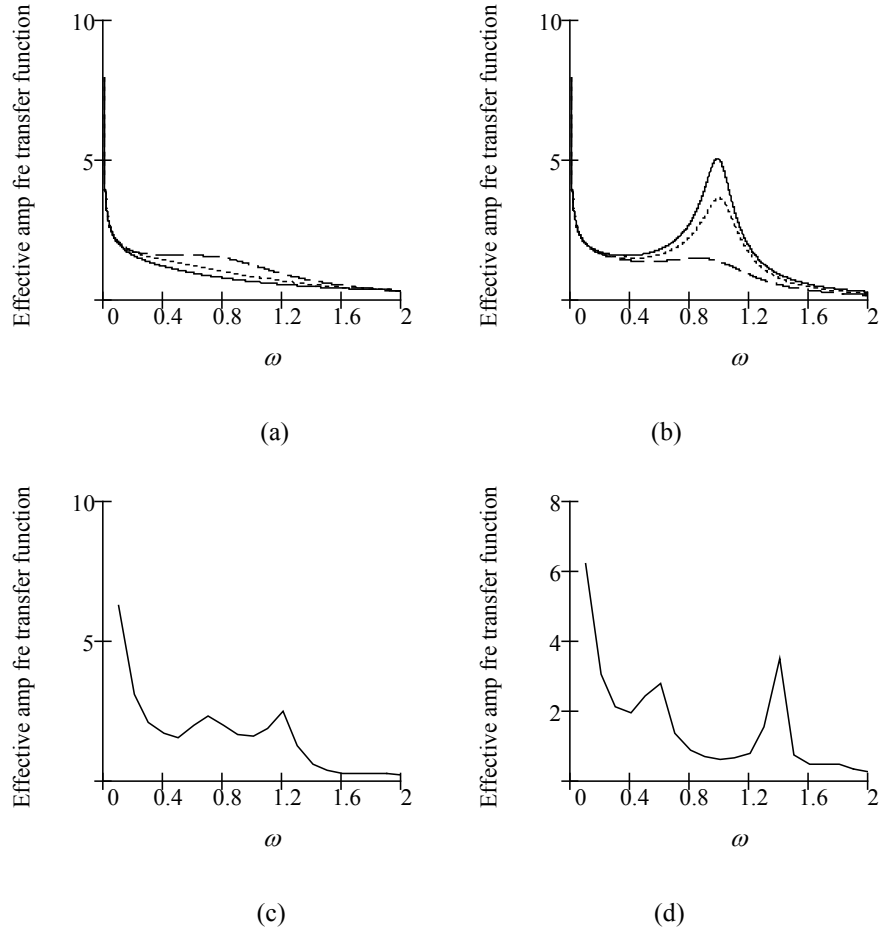

Fig. 5.3. Plots of $|H(\omega)|$ when $m = 1$ and $k = 1$. (a). $|H(\omega)|$ with constant fractional orders for $\omega = 1.1$ and $c = 0.2$, $\beta = 1.3$, $\lambda = 0.3$. Solid: $\alpha = 1.3$. Dot: $\alpha = 1.6$. Dash: $\alpha = 1.9$. (b). $|H(\omega)|$ with constant fractional orders for $\omega = 1.1$ and $c = 0.2$, $\beta = 1.3$, $\lambda = 0.3$. Solid: $\alpha = 2.3$. Dot: $\alpha = 2.6$. Dash: $\alpha = 2.9$. (c). $|H(\omega)|$ with variable fractional orders for $\omega \in (0, 2)$ and $c = 0.2$, $\alpha(\omega) = 1.10 + 1.89|\cos(\omega)|$, $\beta(\omega) = 1 + 0.99|\sin(\omega)|$, $\lambda(\omega) = 0.99|\cos(\omega)|$. (d). $|H(\omega)|$ with variable fractional orders for $\omega \in (0, 2)$ and $c = 1.2$, $\alpha(\omega) = 1.10 + 1.89|\cos(\omega)|$, $\beta(\omega) = 1 + 0.99|\sin(\omega)|$, $\lambda(\omega) = 0.99|\cos(\omega)|$.

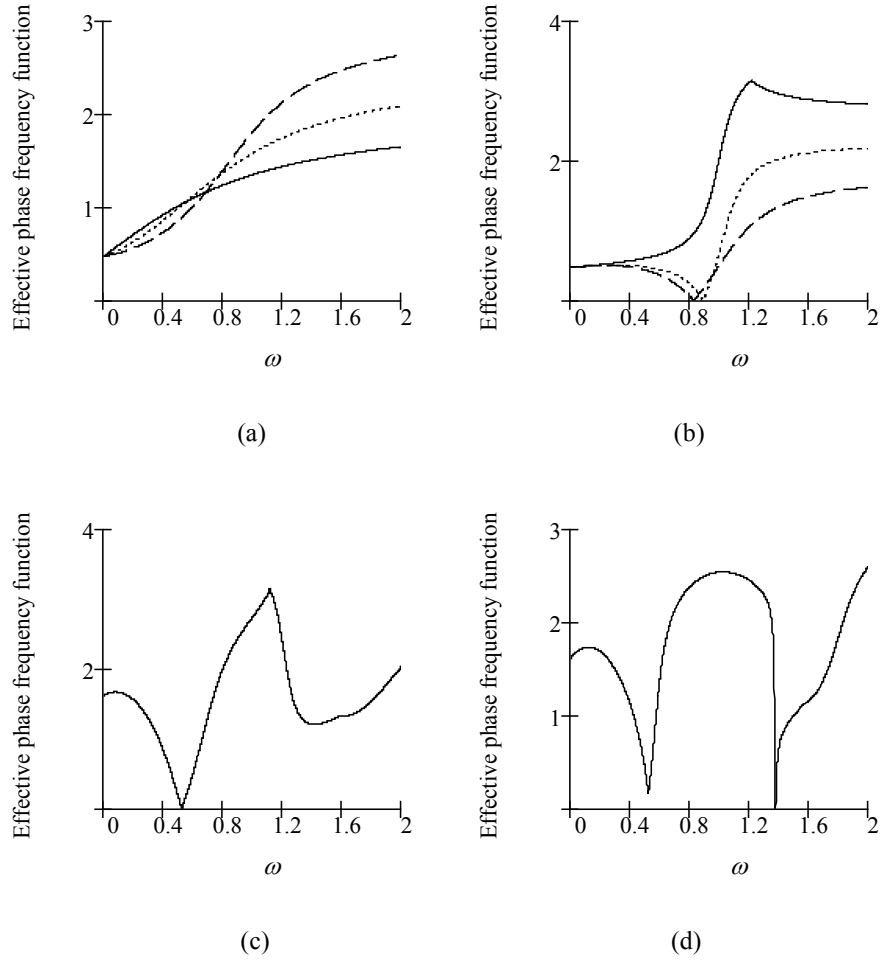

Fig. 5.4. Plots of $\varphi(\omega)$ when $m = 1$, $k = 1$. (a). $\varphi(\omega)$ with constant fractional orders for $\omega = 1.1$ and $c = 0.2$, $\beta = 1.3$, $\lambda = 0.3$. Solid: $\alpha = 1.3$. Dot: $\alpha = 1.6$. Dash: $\alpha = 1.9$. (b). $\varphi(\omega)$ with constant fractional orders for $\omega = 1.1$ and $c = 0.2$, $\beta = 1.3$, $\lambda = 0.3$. Solid: $\alpha = 2.3$. Dot: $\alpha = 2.6$. Dash: $\alpha = 2.9$. (c). $\varphi(\omega)$ with variable fractional orders for $\omega \in (0, 2)$ and $c = 0.2$, $\alpha(\omega) = 1.10 + 1.89|\cos(\omega)|$, $\beta(\omega) = 1 + 0.99|\sin(\omega)|$, $\lambda(\omega) = 0.99|\cos(\omega)|$. (d). $\varphi(\omega)$ with variable fractional orders for $\omega \in (0, 2)$ and $c = 1.2$, $\alpha(\omega) = 1.10 + 1.89|\cos(\omega)|$, $\beta(\omega) = 1 + 0.99|\sin(\omega)|$, $\lambda(\omega) = 0.99|\cos(\omega)|$.

## 6. Mathematical explanation of Rayleigh damping assumption

Rayleigh introduced his damping assumption in [92]. It has been widely used in structural mechanics, see e.g., Palley et al. [90], Trombetti and Silvestri [93], Poul and Zerva [94], Nåvik et al. [95], Park et al. [96], Hussein et al. [97], Sigmund and Jensen [98], Chen et al. [99], Battisti et al. [100], Cox et al. [101], Tisseur and Meerbergen [102], Chu et al. [103], Fay et al. [104], Naderian et al. [105], Tian et al. [106], Iovane et al. [107], Kouris et al. [108], Jin and Xia [109], just mentioning a few. In the field, it is well known that how to give a mathematical explanation of the Rayleigh damping assumption is an open problem.

Recently, the author presented a mathematical explanation of the Rayleigh damping assumption when $\alpha$, $\beta$, and $\lambda$ are constants [24, Chap. 14]. In this paper, we address the further explanation by taking into account the functions $\alpha(\omega)$, $\beta(\omega)$, and $\lambda(\omega)$ instead of $\alpha$, $\beta$, and $\lambda$ being constants.

Let $c_r$ be the Rayleigh damping. Its standard form is expressed by

$$c_{\text{ray}} = am + bk, \tag{6.1}$$

where $a$ and $b$ are frequency-dependent parameters. Precisely, $a$ is proportional to $\omega$ while $b$ is inversely proportional to $\omega$. With the variable-order fractional vibrator (1.12), we use its effective damping to give a more general explanation of the Rayleigh damping assumption.

**Theorem 6.1.** Let $c = 0$ in (2.11). Then, $c_{\text{eff}}$ in (2.11) reduces to the generalized Rayleigh damping, denoted by $c_{\text{gray}}$, in the form

$$c_{\text{gray}} = m\omega^{\alpha(\omega)-1}\sin\frac{\alpha(\omega)\pi}{2} + k\omega^{\lambda(\omega)-1}\sin\frac{\lambda(\omega)\pi}{2}. \tag{6.2}$$

*Proof.* Let $c = 0$ in (2.11). Then, $c_{\text{eff}}$ becomes the above. The above exhibits that $c_{\text{gray}}$ is proportional to $m$ with the coefficient $a$ in the form

$$a = a(\omega) = \omega^{\alpha(\omega)-1}\sin\frac{\alpha(\omega)\pi}{2}. \tag{6.3}$$

Besides, $c_{\text{gray}}$ is proportional to $k$ with the coefficient $b$ in the form

$$b = b(\omega) = \omega^{\lambda(\omega)-1}\sin\frac{\lambda(\omega)\pi}{2}. \tag{6.4}$$

Thus, we have

$$c_{\text{gray}} = a(\omega)m + b(\omega)k. \tag{6.5}$$

Because $a(\omega)$ is proportional to $\omega$ while $b(\omega)$ is inversely proportional to $\omega$. Eq. (6.2) is a mathematical expression that is consistent with the Rayleigh damping assumption. □

Note that $a(\omega) \geq 0$ and $b(\omega) \geq 0$ for $1 < \alpha(\omega) < 2$ and $0 \leq \lambda(\omega) < 1$. Accordingly, $c_{\text{gray}} \geq 0$ in that cases of $1 < \alpha(\omega) < 2$ and $0 \leq \lambda(\omega) < 1$. On the other hand, if $2 < \alpha(\omega) < 3$ and $0 \leq \lambda(\omega) < 1$, $a(\omega)$ may be negative and $b(\omega) \geq 0$. In that case, $c_{\text{gray}}$ may be negative. Because $\alpha(\omega)$ and $\lambda(\omega)$ are varying in terms of $\omega$, (6.2) is a general form of the Rayleigh damping. In fact, when both $\alpha(\omega)$ and $\lambda(\omega)$ are constants $\alpha$ and $\lambda$, $c_{\text{gray}}$ reduces to a specific form of the standard Rayleigh damping in the form

$$c_r = m\omega^{\alpha-1}\sin\frac{\alpha\pi}{2} + k\omega^{\lambda-1}\sin\frac{\lambda\pi}{2}. \tag{6.6}$$

Fig. 6.1 shows the plots of $a(\omega)$ and $b(\omega)$ while Fig. 6.2 demonstrates the plots of $c_{\text{gray}}$.

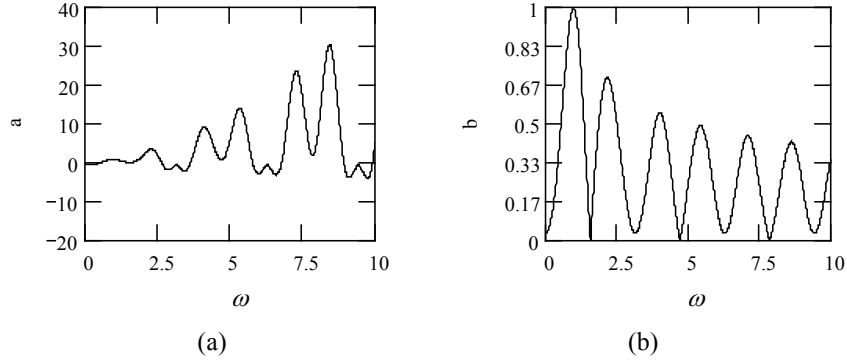

Fig. 6.1. Plots of $a(\omega)$ and $b(\omega)$ with $m = k = 1$. (a). $a(\omega)$ with variable fractional order. $\alpha(\omega) = 1.80 + 1.19|\sin(\omega)|$. (b). $b(\omega)$ with variable fractional order. $\lambda(\omega) = 0.99|\cos(\omega)|$.

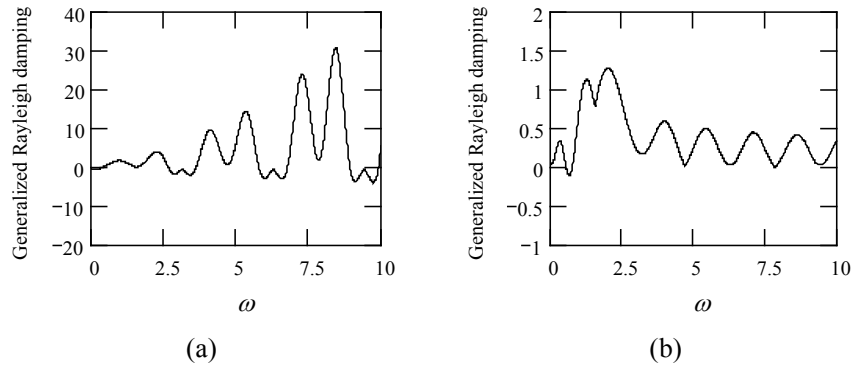

Fig. 6.2. Plots of $c_{\text{gray}}$ with $m = k = 1$. (a). $c_{\text{gray}}$ with variable fractional orders. $\alpha(\omega) = 1.80 + 1.19|\sin(\omega)|$ and $\lambda(\omega) = 0.99|\cos(\omega)|$. (b). $c_{\text{gray}}$ with variable fractional orders. $\alpha(\omega) = 2.99\exp(-\omega)$ and $\lambda(\omega) = 0.99|\cos(\omega)|$.

## 7. Results for the variable-order fractional vibrators (1.7) – (1.11)

We now propose the results regarding the variable-order fractional vibrators expressed by (1.7) – (1.11). Mathematically, (1.7) – (1.11) are the special cases of (1.12). However, from a view of engineering, studying each class from (1.7) to (1.11) is meaningful since each has its specific application area. For instance, (1.7) is for the case of damping free structures with variable-order fractional inertia force and conventional restoration one. The class (1.10) stands for the case of damping free structures with variable-order fractional inertia force and variable-order fractional restoration one. The importance of the results in Section 7 lies in representing the results for each class from (1.7) – (1.11) based on the results with respect to (1.12) addressed in Sections 2 – 5. In what follows, $x_0$ and $v_0$ are initial conditions.

### 7.1. Results regarding variable-order fractional vibrator (1.7)

Consider (1.7). Letting $c = 0$ and $\lambda(\omega) = 0$ in (2.2) yields the equivalent equation of (1.7) in the form

$$-m\omega^{\alpha(\omega)-2}\cos\frac{\alpha(\omega)\pi}{2}\frac{d^2x(t)}{dt^2}+m\omega^{\alpha(\omega)-1}\sin\frac{\alpha(\omega)\pi}{2}\frac{dx(t)}{dt}+kx(t)=f(t),\ \ 1<\alpha(\omega)<3. \tag{7.1}$$

Let $m_{\text{eff1}}$ be the effective mass of the fractional vibrator (1.7). Then, letting $c = 0$ in (2.10) produces

$$m_{\text{eff1}} = -m\omega^{\alpha(\omega)-2}\cos\frac{\alpha(\omega)\pi}{2}. \tag{7.2}$$

Denote by $c_{\text{eff1}}$ the effective damping of the fractional vibrator (1.7). Then, letting $c = \lambda(\omega) = 0$ in (2.11) results in

$$c_{\text{eff1}} = m\omega^{\alpha(\omega)-1}\sin\frac{\alpha(\omega)\pi}{2}. \tag{7.3}$$

Let $\omega_{\text{effn1}}$ be the restricted effective damping free natural frequency of (1.7) for $m_{\text{eff1}} > 0$. It is given by

$$\omega_{\text{effn1}} = \sqrt{\frac{k}{m_{\text{eff1}}}} = \sqrt{\frac{k}{-m\omega^{\alpha(\omega)-2}\cos\frac{\alpha(\omega)\pi}{2}}}. \tag{7.4}$$

Let $\zeta_{\text{eff1}}$ be the restricted effective damping ratio of (1.7) under the condition of $m_{\text{eff1}} > 0$. Then, from (4.1), we have

$$\zeta_{\text{eff1}} = \frac{c_{\text{eff1}}}{2\sqrt{m_{\text{eff1}}k}} = \frac{m\omega^{\alpha(\omega)-1}\sin\frac{\alpha(\omega)\pi}{2}}{2\sqrt{-mk\omega^{\alpha(\omega)-2}\cos\frac{\alpha(\omega)\pi}{2}}}. \tag{7.5}$$

Consider the restricted effective damped natural frequency for (1.7) when $m_{\text{eff1}} > 0$. Denote it by $\omega_{\text{effd1}}$. With the restriction $|\zeta_{\text{eff1}}| \le 1$, we have

$$\omega_{\text{effd1}} = \omega_{\text{effn1}}\sqrt{1-\zeta_{\text{eff1}}^2}. \tag{7.6}$$

Let $\gamma_{\text{eff1}}$ be the restricted effective frequency ratio of (1.7). It is, according to (4.8), given by

$$\gamma_{\text{eff1}} = \frac{\omega}{\omega_{\text{effn1}}}. \tag{7.7}$$

Both $\omega_{\text{effn}}$ and $\gamma_{\text{eff}}$ are restricted by $m_{\text{eff}} > 0$.

Denote by $x_1(t)$ the free response to (1.7). According to Theorem 5.1, we have

$$x_1(t) = e^{-\zeta_{\text{eff1}}\omega_{\text{effn1}}t}\left(x_0\cos\omega_{\text{effd1}}t + \frac{v_0+\zeta_{\text{eff1}}\omega_{\text{effn1}}x_0}{\omega_{\text{effd1}}}\sin\omega_{\text{effd1}}t\right),\ \ t \ge 0. \tag{7.8}$$

Let $h_1(t)$ be the impulse response to (1.7). Based on Theorem 5.2, we have

$$h_1(t) = e^{-\varsigma_{\text{eff1}}\omega_{\text{effn1}}t} \frac{1}{m_{\text{eff1}}\omega_{\text{effd1}}} \sin\omega_{\text{effd1}}t, \quad t \geq 0. \tag{7.9}$$

Let $H_1(\omega)$ be the restricted frequency transfer function of the fractional vibrator (1.7). Following Theorem 5.3, we have

$$H_1(\omega) = \frac{1}{k\left(1 - \gamma_{\text{eff1}}^2 + i2\varsigma_{\text{eff1}}\gamma_{\text{eff1}}\right)}. \tag{7.10}$$

**7.2. Results regarding variable-order fractional vibrator (1.8)**

Taking into account (1.8), we substitute $\alpha(\omega)$ by 2 and $\lambda(\omega)$ with 0 in (2.2). Then, we have the equivalent equation of (1.8) given by

$$\left[m - c\omega^{\beta(\omega)-2}\cos\frac{\beta(\omega)\pi}{2}\right]\frac{d^2x(t)}{dt^2} + c\omega^{\beta(\omega)-1}\sin\frac{\beta(\omega)\pi}{2}\frac{dx(t)}{dt} + kx(t) = f(t), \quad 0 < \beta(\omega) < 2. \tag{7.11}$$

Write the above by

$$m_{\text{eff2}}\frac{d^2x(t)}{dt^2} + c_{\text{eff2}}\frac{dx(t)}{dt} + kx(t) = f(t), \tag{7.12}$$

where $m_{\text{eff2}}$ and $c_{\text{eff2}}$ are the effective mass and damping of the fractional vibrator (1.8), respectively. From (7.11) and (7.12), we see that

$$m_{\text{eff2}} = m - c\omega^{\beta(\omega)-2}\cos\frac{\beta(\omega)\pi}{2}, \tag{7.13}$$

$$c_{\text{eff2}} = c\omega^{\beta(\omega)-1}\sin\frac{\beta(\omega)\pi}{2}. \tag{7.14}$$

Denote by $\omega_{\text{effn2}}$ the restricted effective damping free natural frequency of (1.8) when $m_{\text{eff2}} > 0$. It is

$$\omega_{\text{effn2}} = \sqrt{\frac{k}{m_{\text{eff2}}}}. \tag{7.15}$$

Let $\varsigma_{\text{eff2}}$ be the restricted effective damping ratio of (1.8) in the case of $m_{\text{eff2}} > 0$. Then,

$$\varsigma_{\text{eff2}} = \frac{c_{\text{eff2}}}{2\sqrt{m_{\text{eff2}}k}}. \tag{7.16}$$

Denote by $\omega_{\text{effd2}}$ the restricted effective damped natural frequency of (1.8) for $m_{\text{eff2}} > 0$. Using the restriction $|\varsigma_{\text{eff2}}| \leq 1$, we have

$$\omega_{\text{effd2}} = \omega_{\text{effn2}}\sqrt{1 - \varsigma_{\text{eff2}}^2}. \tag{7.17}$$

Let $\gamma_{\text{eff2}}$ be the restricted effective frequency ratio of (1.8). It is given by

$$\gamma_{\text{eff2}} = \frac{\omega}{\omega_{\text{effn2}}}. \tag{7.18}$$

Let $x_2(t)$ be the free response to (1.8). With Theorem 5.1 or from (7.12), we have

$$x_2(t) = e^{-\varsigma_{\text{eff2}}\omega_{\text{effn2}}t}\left(x_0 \cos\omega_{\text{effd2}}t + \frac{v_0 + \varsigma_{\text{eff2}}\omega_{\text{effn2}}x_0}{\omega_{\text{effd2}}}\sin\omega_{\text{effd2}}t\right), \quad t \geq 0. \tag{7.19}$$

Denote by $h_2(t)$ the impulse response to (1.8). Following Theorem 5.2 or from (7.12), we have

$$h_2(t) = e^{-\varsigma_{\text{eff2}}\omega_{\text{effn2}}t}\frac{1}{m_{\text{eff2}}\omega_{\text{effd2}}}\sin\omega_{\text{effd2}}t, \quad t \geq 0. \tag{7.20}$$

Let $H_2(\omega)$ be the restricted frequency transfer function of the fractional vibrator (1.8). With Theorem 5.3 or by doing the Fourier transform of $h_2(t)$, we have

$$H_2(\omega) = \frac{1}{k\left(1 - \gamma_{\text{eff2}}^2 + i2\varsigma_{\text{eff2}}\gamma_{\text{eff2}}\right)}. \tag{7.21}$$

**7.3. Results regarding variable-order fractional vibrator (1.9)**

Consider (1.9). Replacing $\lambda(\omega)$ with 0 in (2.2) results in

$$\begin{aligned}&-\left[m\omega^{\alpha(\omega)-2}\cos\frac{\alpha(\omega)\pi}{2} + c\omega^{\beta(\omega)-2}\cos\frac{\beta(\omega)\pi}{2}\right]\frac{d^2x(t)}{dt^2}\\&+\left[m\omega^{\alpha(\omega)-1}\sin\frac{\alpha(\omega)\pi}{2} + c\omega^{\beta(\omega)-1}\sin\frac{\beta(\omega)\pi}{2}\right]\frac{dx(t)}{dt}\\&+kx(t) = f(t), \quad 1 < \alpha(\omega) < 3, \ 0 < \beta(\omega) < 2.\end{aligned} \tag{7.22}$$

The above is the equivalent equation of (1.9). Let $m_{\text{eff3}}$ and $c_{\text{eff3}}$ be the effective mass and damping of the fractional vibrator (1.9), respectively. Write the above by

$$m_{\text{eff3}}\frac{d^2x(t)}{dt^2} + c_{\text{eff3}}\frac{dx(t)}{dt} + kx(t) = f(t), \tag{7.23}$$

By comparing (7.22) to (7.23), we have

$$m_{\text{eff3}} = -\left[m\omega^{\alpha(\omega)-2}\cos\frac{\alpha(\omega)\pi}{2} + c\omega^{\beta(\omega)-2}\cos\frac{\beta(\omega)\pi}{2}\right], \tag{7.24}$$

$$c_{\text{eff3}} = m\omega^{\alpha(\omega)-1}\sin\frac{\alpha(\omega)\pi}{2} + c\omega^{\beta(\omega)-1}\sin\frac{\beta(\omega)\pi}{2}. \tag{7.25}$$

Let $\omega_{\text{effn3}}$ be the restricted effective damping free natural frequency of (1.9) for $m_{\text{eff3}} > 0$. It is given by

$$\omega_{\text{effn3}} = \sqrt{\frac{k}{m_{\text{eff3}}}}. \tag{7.26}$$

Let $\zeta_{\text{eff3}}$ be the restricted effective damping ratio of (1.9) when $m_{\text{eff3}} > 0$. Then,

$$\zeta_{\text{eff3}} = \frac{c_{\text{eff3}}}{2\sqrt{m_{\text{eff3}} k}}. \tag{7.27}$$

Denote by $\omega_{\text{effd3}}$ the restricted effective damped natural frequency of (1.9) when $m_{\text{eff3}} > 0$. With $|\zeta_{\text{eff3}}| \leq 1$, we have

$$\omega_{\text{effd3}} = \omega_{\text{effn3}} \sqrt{1 - \zeta_{\text{eff3}}^2}. \tag{7.28}$$

Denote by $\gamma_{\text{eff3}}$ the restricted effective frequency ratio of (1.9). It is in the form

$$\gamma_{\text{eff3}} = \frac{\omega}{\omega_{\text{effn3}}}. \tag{7.29}$$

Denote $x_3(t)$ the free response to (1.9). Since (7.23) is a standard vibration equation, $x_3(t)$ is given by

$$x_3(t) = e^{-\zeta_{\text{eff3}} \omega_{\text{effn3}} t} \left( x_0 \cos \omega_{\text{effd3}} t + \frac{v_0 + \zeta_{\text{eff3}} \omega_{\text{effn3}} x_0}{\omega_{\text{effd3}}} \sin \omega_{\text{effd3}} t \right), \quad t \geq 0. \tag{7.30}$$

Denote by $h_3(t)$ the impulse response to (1.9). Following Theorem 5.2 or from (7.23), we have

$$h_3(t) = e^{-\zeta_{\text{eff3}} \omega_{\text{effn3}} t} \frac{1}{m_{\text{eff3}} \omega_{\text{effd3}}} \sin \omega_{\text{effd3}} t, \quad t \geq 0. \tag{7.31}$$

Let $H_3(\omega)$ be the restricted frequency transfer function of the fractional vibrator (1.9). Doing the Fourier transform of $h_3(t)$ yields

$$H_3(\omega) = \frac{1}{k \left( 1 - \gamma_{\text{eff3}}^2 + i 2 \zeta_{\text{eff3}} \gamma_{\text{eff3}} \right)}. \tag{7.32}$$

**7.4. Results regarding variable-order fractional vibrator (1.10)**

Let the equivalent equation of (1.10) be

$$m_{\text{eff4}} \frac{d^2 x(t)}{dt^2} + c_{\text{eff4}} \frac{dx(t)}{dt} + k_{\text{eff4}} x(t) = f(t), \tag{7.33}$$

where $m_{\text{eff4}}$, $c_{\text{eff4}}$, and $k_{\text{eff4}}$ are the effective mass, damping, and stiffness of the fractional vibrator (1.10), respectively. Taking into account (1.10), we replace $c$ by 0 in (2.2) and have the equivalent equation of (1.10) expressed by

$$-m\omega^{\alpha(\omega)-2}\cos\frac{\alpha(\omega)\pi}{2}\frac{d^2x(t)}{dt^2}+\left[m\omega^{\alpha(\omega)-1}\sin\frac{\alpha(\omega)\pi}{2}+k\omega^{\lambda(\omega)-1}\sin\frac{\lambda(\omega)\pi}{2}\right]\frac{dx(t)}{dt}$$
$$+k\omega^{\lambda(\omega)}\cos\frac{\lambda(\omega)\pi}{2}x(t)=f(t),\quad 1<\alpha(\omega)<3,\ 0\le\lambda(\omega)<1. \tag{7.34}$$

By comparing (7.33) to (7.34), we see that

$$m_{\text{eff4}}=-m\omega^{\alpha(\omega)-2}\cos\frac{\alpha(\omega)\pi}{2}, \tag{7.35}$$

$$c_{\text{eff4}}=m\omega^{\alpha(\omega)-1}\sin\frac{\alpha(\omega)\pi}{2}+k\omega^{\lambda(\omega)-1}\sin\frac{\lambda(\omega)\pi}{2}, \tag{7.36}$$

$$k_{\text{eff4}}=k\omega^{\lambda(\omega)}\cos\frac{\lambda(\omega)\pi}{2}. \tag{7.37}$$

Because (7.33) is a standard vibration equation in form, we write its free response, denoted by $x_4(t)$, by

$$x_4(t)=e^{-\varsigma_{\text{eff4}}\omega_{\text{effn4}}t}\left(x_0\cos\omega_{\text{effd4}}t+\frac{v_0+\varsigma_{\text{eff4}}\omega_{\text{effn4}}x_0}{\omega_{\text{effd4}}}\sin\omega_{\text{effd4}}t\right),\quad t\ge 0, \tag{7.38}$$

where $\varsigma_{\text{eff4}}=\dfrac{c_{\text{eff4}}}{2\sqrt{m_{\text{eff4}}k_{\text{eff4}}}}$ is the restricted effective damping ratio of (1.10), $\omega_{\text{effn4}}=\sqrt{\dfrac{k_{\text{eff4}}}{m_{\text{eff4}}}}$ is the restricted effective damping free natural frequency of (1.10), $\omega_{\text{effd4}}=\omega_{\text{effn4}}\sqrt{1-\varsigma_{\text{eff4}}^2}$ for $|\varsigma_{\text{eff4}}|\le 1$ is the restricted effective damped natural frequency of (1.10). By restricted, we mean that they are in the sense of $m_{\text{eff4}}>0$.

Denote by $h_4(t)$ the impulse response to (1.10). From Theorem 5.2 or (7.33), we at once write it by

$$h_4(t)=e^{-\varsigma_{\text{eff4}}\omega_{\text{effn4}}t}\frac{1}{m_{\text{eff4}}\omega_{\text{effd4}}}\sin\omega_{\text{effd4}}t,\quad t\ge 0. \tag{7.39}$$

According to Theorem 5.3, the restricted frequency transfer function of the fractional vibrator (1.10), which is denoted by $H_4(\omega)$, is given by

$$H_4(\omega)=\frac{1}{k_{\text{eff4}}\left(1-\gamma_{\text{eff4}}^2+i2\varsigma_{\text{eff4}}\gamma_{\text{eff4}}\right)}, \tag{7.40}$$

where $\gamma_{\text{eff4}}=\dfrac{\omega}{\omega_{\text{effn4}}}$.

**7.5. Results regarding variable-order fractional vibrator (1.11)**

The equivalent equation of (1.11) is expressed by

$$m\frac{d^2x(t)}{dt^2} + k\omega^{\lambda(\omega)-1}\sin\frac{\lambda(\omega)\pi}{2}\frac{dx(t)}{dt} + k\omega^{\lambda(\omega)}\cos\frac{\lambda(\omega)\pi}{2}x(t) = f(t), \quad 0 \le \lambda(\omega) < 1. \tag{7.41}$$

As a matter of fact, in (2.2), letting $\alpha(\omega) = 2$ and $c = 0$ produces the above.

Let $c_{\text{eff5}}$ and $k_{\text{eff5}}$ be the effective damping and stiffness of the fractional vibrator (1.11), respectively. From the above, we have

$$c_{\text{eff5}} = k\omega^{\lambda(\omega)-1}\sin\frac{\lambda(\omega)\pi}{2}, \tag{7.42}$$

$$k_{\text{eff5}} = k\omega^{\lambda(\omega)}\cos\frac{\lambda(\omega)\pi}{2}. \tag{7.43}$$

Thus, (7.41) can be rewritten by

$$m\frac{d^2x(t)}{dt^2} + c_{\text{eff5}}\frac{dx(t)}{dt} + k_{\text{eff5}}x(t) = f(t). \tag{7.44}$$

The above in form is a standard vibration equation. Let $x_5(t)$ be the free response of (7.44). Then,

$$x_5(t) = e^{-\varsigma_{\text{eff5}}\omega_{\text{effn5}}t}\left(x_0\cos\omega_{\text{effd5}}t + \frac{v_0 + \varsigma_{\text{eff5}}\omega_{\text{effn5}}x_0}{\omega_{\text{effd5}}}\sin\omega_{\text{effd5}}t\right), \quad t \ge 0, \tag{7.45}$$

where $\varsigma_{\text{eff5}} = \dfrac{c_{\text{eff5}}}{2\sqrt{mk_{\text{eff5}}}}$, $\omega_{\text{effn5}} = \sqrt{\dfrac{k_{\text{eff5}}}{m}}$, and $\omega_{\text{effd5}} = \omega_{\text{effn5}}\sqrt{1-\varsigma_{\text{eff5}}^2}$ for $|\varsigma_{\text{eff5}}| \le 1$.

Denote by $h_5(t)$ the impulse response to (1.11). Then,

$$h_5(t) = e^{-\varsigma_{\text{eff5}}\omega_{\text{effn5}}t}\frac{1}{m\omega_{\text{effd5}}}\sin\omega_{\text{effd5}}t, \quad t \ge 0. \tag{7.46}$$

Let $H_5(\omega)$ be the frequency transfer function of the fractional vibrator (1.11). Then,

$$H_5(\omega) = \frac{1}{k_{\text{eff5}}\left(1 - \gamma_{\text{eff5}}^2 + i2\varsigma_{\text{eff5}}\gamma_{\text{eff5}}\right)}, \tag{7.47}$$

where $\gamma_{\text{eff5}} = \dfrac{\omega}{\omega_{\text{effn5}}}$.

**8. Conclusions**

We have derived out an equivalent vibration equation of a class of variable-order fractional vibrators (1.12), which is expressed by (2.2). Base on it, we have presented the analytical expressions of its effective mass, and damping, stiffness in Section 2. We have brought forward the analytical expressions of its damping ratio, damping free natural frequency, damped natural frequency, and frequency ratio in Section 4.

In Section 5, we have put forward the close form expressions of free response, impulse response, and frequency transfer function. The asymptotic properties of its effective mass, damping, and stiffness have been given in Section 3. In addition, we have given a mathematical explanation of the Rayleigh damping assumption when $\alpha(\omega)$ and $\lambda(\omega)$ are varying in terms of $\omega$ in Section 6 Besides, based on (2.2), we have proposed equivalent equations, effective vibration parameters (mass, damping, stiffness, damping ratio, damping free natural frequency, damped natural frequency, frequency ratio), as well as free responses, impulse responses, and frequency transfer functions of five classes of variable-order fractional vibrators (1.7) – (1.11) in Section 7. When $\alpha(\omega)$, $\beta(\omega)$ and $\lambda(\omega)$ are constants, the present theory is applicable to the fractional vibrators expressed by (1.1) – (1.6)[6].

**Declaration of competing interest**

The author declares that there are no known competing financial interests or personal relationships that could have appeared to influence the work reported in this paper.

**Acknowledgements**

This work was supported in part by the National Natural Science Foundation of China (NSFC) under the project grant number 61672238. The views and conclusions contained in this document are those of the author and should not be interpreted as representing the official policies, either expressed or implied, of NSFC or the Chinese government.

---

[6] Refer the author's recent work [24] for the details about (1.1) – (1.6).